\documentclass[11pt, a4paper]{article}
\usepackage{amsmath,amssymb,amsfonts}
\usepackage{a4wide}
\usepackage{epic,eepic,graphicx}
\newtheorem{theorem}{Theorem}
\newtheorem{lemma}[theorem]{Lemma}
\newtheorem{proposition}[theorem]{Proposition}

\newtheorem{definition}{Definition}

\makeatletter

\@addtoreset{equation}{section}
\makeatother
\newcommand{\N}{{\mathbb N}}

\newcommand{\R}{{\mathbb R}}

\newcommand{\pa}{{\partial}}
\newcommand{\na}{{\nabla}}

\newcommand{\eps}{{\varepsilon}}

\def\div{\hbox{div  }}

\def\bmo{{\rm BMO}}
\def\O{{\cal O}}
\renewcommand 		\leq 		{\leqslant}

\title{Regularity issues in the problem of fluid structure interaction}
\author{David G\'erard-Varet
\footnote{DMA/CNRS, Ecole Normale Sup«erieure, 45 rue dÕUlm,  75005 Paris, France}
\and Matthieu Hillairet
\footnote{Universit\'e de Toulouse, IMT \'equipe MIP, 118, route de Narbonne, 31062 Toulouse cedex, France}
}

\begin{document}

\maketitle

\abstract{We investigate the evolution of rigid bodies in a
  viscous incompressible fluid. The flow is governed by the 2D Navier-Stokes
  equations,  set in a bounded domain with Dirichlet boundary
  conditions. The boundaries of the solids and the domain have H\"older
  regularity $C^{1, \alpha}$, $0 < \alpha \le 1$. First, we show the
  existence and uniqueness of strong solutions up to collision. A key
  ingredient is a $\bmo$ bound on the velocity gradient, which substitutes to
  the standard $H^2$ estimate for smoother domains. Then, we
   study the asymptotic behaviour of one $C^{1, \alpha}$ body falling 
   over a flat surface. We show that collision is possible  in finite
  time if and only if $\alpha < 1/2$. 
}

\section{Introduction}
To understand the dynamics of solid bodies immersed in a fluid is  of
primary physical interest, with regards to a wide range of phenomena
such as sedimentation, filtration, or coagulation. For two-dimensional
flows, under the assumption that the $N$  bodies are rigid and homogeneous,  
and that the fluid is  incompressible and viscous, one considers
classically the following model:
\begin{description}
\item[i)]
The velocity  $u$ and pressure $p$ satisfy 
Navier-Stokes equations in the fluid domain $F(t)$:
\begin{equation} \label{NS}
\begin{aligned}
& \overline{\rho} \left( \pa_t u + u \cdot \na u \right) 
 - \mu  \Delta u   = - \na p +  \overline{\rho} f,\\
& \div u =  0,
\end{aligned}
\quad  x \in F(t).
\end{equation}
\item[ii)]
The $N$ solid bodies are described by the closures $\overline{S^i(t)}$
 of connected  bounded domains $S^i(t)$, $1 \le i
\le N$. They have   rigid velocity fields 
\begin{equation} \label{rigidvel}
 u^i(t,x) = v^i(t) + \omega^i(t)  (x - x^i(t))^\bot, \quad
 x \in \overline{S^i(t)}, \quad 1 \le i \le N,
\end{equation}
where $v^i$ and $\omega^i$ are the translation and angular velocities,
whereas $x^i$ is the position of the center of mass.
 \item[iii)]
The moving fluid and solid domains occupy a fixed bounded domain
$\Omega$ of $\R^2$,  with Dirichlet boundary condition:  
\begin{equation} \label{dirichlet}
F(t) = \Omega \setminus \cup_{i=1}^N \overline{S^i(t)}, \quad 
  u = 0, \quad x \in \pa \Omega.
\end{equation}
\item[iv)]
The fluid and solid systems are coupled by  the  continuity of the velocity, 
\begin{equation} \label{contin1}
u = u^i, \quad x \in \pa S^i(t),
\end{equation}
and the continuity of the stresses: 
\begin{equation} \label{contin2}
\begin{aligned}
& m^i \, \dot{v}^i(t) = \int_{\pa S^i(t)} \left( \mu \frac{\pa u}{\pa n } - p
  n\right) \, d\sigma \: + \: \int_{S^i(t)} \overline{\rho}^i f, \\
&  J^i \, \dot{\omega}^i(t) = \int_{\pa S^i(t)} (x - x^i)^\bot \cdot 
\left(\mu \frac{\pa u}{\pa n } - p n\right) \, d\sigma \: + \:
\int_{S^i(t)} (x - x^i)^\bot \cdot  \overline{\rho}^i f.
\end{aligned}
\end{equation}
\end{description}
The positive constants $\overline{\rho}$,  $\mu$ are the density and
viscosity of the fluid. The positive constants $m^i$, $J^i$, $\overline{\rho}^i$
are  the total  mass, moment of inertia and density of the $i$-th
solid. The source term $f$ models an  additional forcing (like
gravity). The vector $n$ at the boundary $\pa U$ of an open set $U$  refers
as usual to the outward unit normal vector. 

\medskip
Although natural, these equations exhibit some unexpected
features,  in both two and three dimensions. 
Hence, consider the  case of one rigid body falling in a
cavity ($N=1$), under the action of gravity. It can be shown that 
 {\em if the boundaries of the body and the cavity
are smooth, then no collision can occur in finite time}. In other
words, this  system predicts that  the kinetic energy of the body 
is strongly dissipated by the viscosity, resulting in no collision between the
body and the boundary. This fact has been known from physicists
 for many years  \cite{Brenner&Cox63,Cooley&Oneill69,Dean&Oneill64}, and was
recently proved in one \cite{Vazquez&Zuazua06} and two dimensions
\cite{Hesla05,Hillairet07}.

This no-collision result is of course paradoxical. At the level of
medium-sized objects, it goes against Archimede's law, and is clearly
denied by common experiments. At a microscopic level, it also lacks
 relevance, as rebounds between particles are often involved. Many
physics papers have been  devoted to this paradox, trying to identify the
flaw of the previous modelling. We refer to  the articles \cite{Barnocky,Davis3} among many. 
Among the possible explanations, one of the most popular is {\em
  roughness}. Indeed, the no-collision result relies
on the fact that the boundary of the solid structure is regular enough
(namely $C^{1,1}$). Small irregularities could then explain the occurence
of collisions, see \cite{Maury,Smart} . Moreover, the effect of surface roughness
in the dynamics of particles has been recently emphasized in
experiments \cite{Davis1,Joseph,Davis2}. 

{\em The aim of this paper is to study mathematically the roughness-induced
effect on the collision process.} Therefore, we  consider
H\"older boundaries. Namely, we assume that  
\begin{equation} 
 \pa \Omega \: \in \:  C^{1, \alpha}, \quad  \pa S^i  \in   C^{1,
   \alpha}, \quad \forall \, i, \quad 0 < \alpha \le 1. 
\end{equation}
We will first consider the well-posedness of system
\eqref{NS}-\eqref{contin2}, for such  boundaries. We will establish
existence and uniqueness of some strong solutions, up to
collision. Our result extends previous results obtained for $C^{1,1}$
boundaries. Once this well-posedness is obtained, we will turn to the
question of collision in finite time. We will consider the special
case of one  $C^{1,\alpha}$ rigid body, falling vertically over a
horizontal flat surface. Losely, we will show the following:
\begin{enumerate}
\item {\em For $\alpha \ge 1/2$, no collision can occur, and the
    strong solution exists for all time}. 
\item {\em For $\alpha < 1/2$, one can find  solutions for which
    collision occurs}.
\end{enumerate}
This sharp criteria illustrates that roughness might be the reason for
collision in fluid structure interaction, and the reason for the
apparent paradox of the classical modelling.

\medskip
Before stating precisely the results, let us mention former
mathematical studies. Fluid-solid interaction has been the subject of
many papers, mostly devoted to the existence theory for problem
\eqref{NS}-\eqref{contin2}. A key ingredient in many existence results
is  a weak  formulation of the equations. Introducing the global 
quantities 
\begin{equation} \label{globalv}
v(t,x) \:  := \:  u(t,x) \,  \mathbf{1}_{F(t)}(x) \, + \, 
 \sum_{i=1}^N u^i(t,x) \, \mathbf{1}_{S^i(t)}(x), 
\end{equation}
\begin{equation} \label{globalrho} 
 \rho(t,x) \:   := \:  \rho^F(t,x) \: + \: \sum_{i=1}^N  \rho^i(t,x) \: := \: 
\overline{\rho} \, \mathbf{1}_{F(t)}(x) \, + \,  
\sum_{i=1}^N \,   \overline{\rho}^i   \mathbf{1}_{S^i(t)}(x), 
\end{equation}
the conservations of global momentum, global mass, and bodies masses 
 yield respectively: 
 for all $T > 0$, for all $\displaystyle 
\varphi \in {\cal V}$, for all $\displaystyle \psi \in
 {\cal D}([0,T) \times \Omega)$,  
\begin{equation} \label{var}
\begin{aligned}
&  \int_0^T \int_\Omega \Bigl( \rho v \cdot \pa_t \varphi  + \rho v
  \otimes  v  :
D(\varphi)  - 2 \mu D(v) : D(\varphi) + \rho f \cdot \varphi  \Bigr) dx
ds \: + \:     \int_\Omega \rho_0 v_0 \cdot 
\varphi(0) \: = \: 0, \\
& \int_0^T \int_{\Omega} \Bigl( \rho \pa_t \psi + \rho u \cdot  \na
\psi \Bigr)  + \int_\Omega \rho_0 \psi(0) \: = \: 0, \\
& \int_0^T \int_{\Omega} \Bigl( \rho^i \pa_t \psi + \rho^i u \cdot  \na
\psi \Bigr)  + \int_\Omega \rho^i_0 \psi(0) \: = \: 0. \\
\end{aligned} 
\end{equation} 
The space of  test functions ${\cal V}$ is
\begin{equation*}
{\cal  V} \: = \: \Bigl\{  \varphi \in {\cal D}([0,T) \times \Omega), \quad 
 \na \cdot  \varphi  = 0, \quad \rho^i(t) D(\varphi) = 0, \: \forall \,
 t, \: \forall \, 1 \le i \le N 
\Bigr\}. 
\end{equation*}
The divergence, rigidity inside the fluid, and no-slip condition read
 respectively: 
\begin{equation} \label{var2}
\na \cdot v =0, \quad  \, \rho^i \, D(v) = 0,  \: 1 \le i \le N,
\quad 
v\vert_{\pa \Omega} = 0. 
\end{equation}
We refer to B. Desjardins and M. Esteban \cite{Desjardins&Esteban99} for
the derivation of these 
equations. Similarly to $\rho$, $\rho^i$  and $v$, the initial data
$\rho_0$, $\rho^i_0$  and
 $v_0$ are built upon the initial positions of the
bodies $S^i_0$ and  the initial fluid and solid velocities $u_0$,
$v^i_0$, $\omega^i_0$. We will  assume that there is
no-contact initially, which means 
\begin{equation} \label{nocontact}
\overline{S^i_0} \cap \overline{S^j_0} = \emptyset,
 \quad \overline{S^i_0} \subset \Omega, \quad \forall \, 1 \le i,j \le N,
 \quad i \neq j.
\end{equation}

Broadly speaking, previous studies deal with two kinds of solutions:
weak and strong. 
\begin{definition} 
A {\em weak solution} on $(0,T)$, $T >0$,  is a family  
$$(S^i(t), F(t), 
v), \quad  1 \le i \le N, \quad 
F(t) = \Omega \setminus \cup_{i=1}^N \overline{S^i(t)}$$ 
such that 
\begin{description}
\item[i)] $S^i(t)$ is a connected  bounded domain, for all $0 < t <
  T$, for all $1 \le i \le N$.
\item[ii)] The scalar functions $\rho$, $\rho^i$ defined in \eqref{globalrho}
  and the vector field $v$ satisfy 
$$ (\rho,\rho^i) \in L^\infty(0, T \times \Omega), \quad v \in L^\infty(0,T; \,
L^2(\Omega)) \cap L^2(0,T; \, H^1_0(\Omega)). $$
and equations \eqref{var}, \eqref{var2}.
\end{description}
\end{definition}
By classical results of R. Di perna and P.-L. Lions \cite{DiPernaLions89} on the 
transport equations (\ref{var}b,c), any $(\rho, \rho^i,v)$ satisfying
ii) has the following additional regularity:
 $$\rho, \rho^i \in C([0,T]; \, L^1(\Omega)),$$
and  the initial data is satisfied in this stronger sense.
 Moreover, any $\rho^i$ satisfying (\ref{var}c) is the characteristic
 function of a measurable set: 
$$\rho^i(t,x) \: = \:  {\bf 1}_{S^i(t)}(x), \quad \mbox{ for a.e. } t,
x$$
see \cite[theorem 2.1, p.23]{Lions96}. However, it is not clear that
 $S^i(t)$ should be open and connected, so that this constraint i) is
 added to the definition of a weak solution. Then, using the rigidity
 condition in  \eqref{var2}, one can
 deduce that $v(t, \cdot)$ is a rigid vector  field on each $S^i(t)$, 
and by  (\ref{var}c),  that  $S^i(t) = {\cal R}_t S^i_0$, for a family
of affine isometries ${\cal R}_t$ Lipschitz in $t$. 

The existence of global in time ($T=+\infty$) weak solutions was
proved by  E. Feireisl \cite{Feireisl03} and San Martin and coauthors
\cite{SanMartin&al02} 
extending earlier studies ``up to collision between solids''
\cite{Desjardins&Esteban00,Hoffmann&Starovoitov99,Hoffmann&Starovoitov00,Conca&al00,Gunzburger&Lee&Seregin00}.   
It  holds in  dimensions 2 and 3, with initial data satisfying
\eqref{nocontact} and  
$$v_0 \in L^2(\Omega), \quad \div
v_0 = 0,  \quad f \in L^2((0,T); \,  H^{-1}(\Omega)). $$
Following  the construction by E. Feireisl, {\em no smoothness of the
  boundaries of the domain and the solids is necessary for the
  existence of weak solutions}. However, the
uniqueness of such solutions is unknown in general,
 even considering  dimension 2 and pre-collisional times. 
 After contact, it is known that uniqueness does not
hold, as some {\em entropy condition} is missing to describe properly
the post-collisional dynamics. This suggests to consider stronger
solutions, namely 
\begin{definition} 
A {\em strong solution} on $(0,T)$ is a weak solution  with the
following additional regularity:  
$$  v \in  L^\infty\left(0,T; H^1_0(\Omega)\right) \cap
L^2\left(0,T; \, W^{1,p}(\Omega)\right) \: \mbox{for all finite} \: p, 
\quad \pa_t v \in L^2(0,T; \, L^2(\Omega)).$$ 
\end{definition}

\medskip
Our first result is the following:
\begin{theorem} {\bf (Well-posedness up to collision)} \label{theo1}

\smallskip 
Let $\displaystyle 
\: v_0 \in H^1_0(\Omega), \quad \rho^i_0 \, D(v_0) = 0, \: \forall \, i, 
 \quad  \: f \in L^2((0,T); \, W^{1,\infty}(\Omega)),
\:  \forall \,T >0.$ 
Assume \eqref{nocontact}, and
$$\displaystyle \pa \Omega \in C^{1, \alpha}, \quad \pa S^i_0 \in C^{1, \alpha},
\quad \forall \, i, \quad 0 < \alpha \le 1.$$
 Then, there exists  a maximal $T_* \in (0, \infty]$ with a unique  strong 
 solution on $(0,T)$ for all $T < T_*$. 
Moreover, this strong solution exists up to the first collision, which
means one of the following alternatives holds true: 
\begin{description} 
\item[i)] $ \displaystyle   T_*=\infty, \quad  \delta(t) > 0 \quad \forall
  \, t.$
\item[ii)] $ \displaystyle  T_* < \infty, \quad   \delta(t) > 0 \quad
  \forall \, t < T_*, \quad 
 \lim_{t \rightarrow T_*} \delta(t) \: = \: 0,$ 
 \end{description}
$$ \mbox{ where } \quad \displaystyle  \delta(t) \: := \: \min \{ \:
d(S^i(t), S^j(t)), \quad  
 d(S^i(t), \pa \Omega), \quad 1 \le i, j \le N, \quad i \neq j\}. $$
 \end{theorem}
Note that by  condition \eqref{nocontact}, and the Lipschitz dependance
of ${\cal R}_t$ described  above, $\delta$ is positive at least for
small times. Our theorem is an extension of results of
B. Desjardins and M. Esteban  \cite{Desjardins&Esteban99}, and  
T. Takahashi \cite{Takahashi03bis}, who proved respectively existence and
uniqueness of strong solutions  in the case $\alpha=1$. 
See also \cite{Grandmont&Maday98} for well-posedness under further
technical assumptions on the solids.
A key argument in these  papers is the classical $L^2 \mapsto  H^2$
regularity  property  for  the inverse of the Stokes operator, which holds  in $C^{1,1}$
 domains. In particular, one can show that   
\begin{equation} \label{ell1} \int_0^T \int_{F(t)} | \na^2 v(t, \cdot) |^2
  < +\infty, \quad 0 < T < 
T^*. 
\end{equation}

In the case of general $C^{1,\alpha}$ domains, this $H^2$ regularity
result is still true away from the boundaries, and $(\rho,v)$ still satisfies
\eqref{NS} in the strong sense, that is for almost every $x,$ $t.$
However Theorem 1 requires a control up to the boundary . We will show that
the following $\bmo$ bound:
$$
\int_0^T  \| \na v(t, \cdot) \|^2_{\bmo(F(t))}  < +\infty, \quad 0 < T <
T^*, 
$$
substitutes to \eqref{ell1}, allowing for our well-posedness result.

\medskip
In a second part, we study if bodies can collide in finite time, that is if
$T_*$ is finite or not. We consider one  $C^{1,\alpha}$ solid that moves
vertically near  a flat 
horizontal surface under the action of gravity. More precisely, let us denote
 $S(t) = {\cal R}_t S_0$ the position of the solid at
time $t$. We make the following assumptions:
\begin{enumerate}
\item The source term  is $f = -ge_2$,  with $g>0$,  $\: e_2 = (0,1)$. 
\item The solid moves along the axis  $x_1 = 0$, that is ${\cal R}_t$
  is a vertical translation.
\item  The only possible collision points are on $x_1 = 0$.  
\item   Near $x_1 = 0$, $\pa \Omega$ is  flat and horizontal
\item Near $x_1 = 0$, the lower and upper parts of $\pa S(t)$ are
  given by
 $$x_2 - x_{-}(t) =  |x_1|^{1+\alpha}, \quad x_2 - x_+(t) = -
 |x_1|^2, \quad 0 < t < T_*.$$ 
\item The solid is heavier than the fluid, \emph{i.e.}, $\rho_{|_{S(t)}}> \rho_{|_{F(t)}}.$
\end{enumerate}
Note that if the initial configuration ($\Omega$, $S_0$, $v_0$) 
is symmetric with respect to the $x_1$-axis,  then the unique strong
solution will be symmetric for all $0<t<T_*$, and the solid 
will  move along the vertical axis. 
Hence, there are plenty of configurations satisfying  1-6. 
A typical one  is shown in figure \ref{fig_typicalsituation}.
Our main result is the following: 

 \begin{figure} \label{fig_typicalsituation}
\begin{center}
\includegraphics[height = 5cm, width=5.5cm]{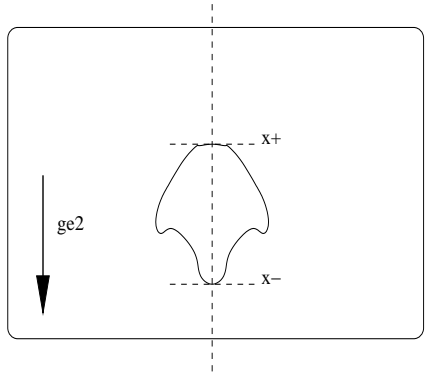}
\end{center}
\caption{Typical situation}
\end{figure}

\begin{theorem} {\bf (Link between collision and boundary
    regularity)} \label{theo2}

\smallskip
For any strong solution satisfying 1-6, $T_*<\infty$ if and only if $\alpha
< 1/2.$
\end{theorem}
In physical terms, the theorem emphasizes the role of roughness in the
collision scenario. Our result extends the results of M. Hillairet 
\cite{Hillairet07} and T.I. Hesla \cite{Hesla05} in the  case 
$\alpha =1$, for which it was shown that no collision occurs. 
Theorem \ref{theo2} relies on the study of the stress $\int_{\pa
  S(t)} (\mu\pa_n u - p \, n)$. When the boundary is regular, this stress
diverges strongly as the distance to the boundary
 goes to zero. This mechanism prevents
collision. When the regularity is weakened, the stress is also
weakened, and  contact may occur.    
The proof of the theorem involves the construction of appropriate test
functions. In that respect, assumptions 2-5 are mostly technical,
allowing  to handle the computations. As can be
seen from our proof, most of our arguments are local, and use only the
{\em weak bounds} given by the conservation of energy. Hence, 
 we believe that, as far as ``real'' (not grazing) collisions are
 concerned, the result might persist for more general domains and weak
 solutions. However, the source term must remain sufficiently integrable,
 as shown by an interesting example of
 Starovoitov \cite{Starovoitov03}. Losely, Starovoitov exhibits an example
 of a weak solution, 
 colliding in finite time, when $\Omega$ and the solid are two
 spheres. But the corresponding source term satisfies only 
$$ f \in L^2(0,T; \, H^{-1}(\Omega)), \quad \forall \, T >0. $$
The $L^2$ norm of $f(t, \cdot)$ diverges as $\delta \rightarrow 0, \: t
\rightarrow T_*$. This allows to compensate the divergence of the
stress and to allow collision, even with regular boundaries. As shown
by the first part of our theorem, this phenomenon is ruled out for more
realistic forcing (such as gravity). 

\medskip
The rest of the article is organized in three sections. 
Section 2 gathers regularity properties for 
the Stokes operator in $C^{1, \alpha}$ domains. Section 3 is devoted
to the proof of Theorem \ref{theo1}. Section 4 contains the proof of
Theorem \ref{theo2}.

\section{Regularity properties in $C^{1, \alpha}$ domains}
Existence and uniqueness of  strong solutions have only
been considered when solids have $C^{1,1}$ boundaries. More precisely, a key
argument in the papers of B. Desjardins and M. Esteban or T. Takahashi
  is  the regularity estimate
\begin{equation} \label{Sobolev}
\| \na u \|_{H^1({\cal O})} + \| p  \|_{H^1({\cal O})/\R} \: \le C \, \left( \| F
\|_{H^1({\cal O})} \: + \:   \| g \|_{H^1({\cal O})} \right) 
\end{equation}
satisfied by the weak solution $(u,p)$ of the Stokes system
\begin{equation} \label{Stokes}
\left\{
\begin{aligned}
 -\Delta u + \na p  & \: =  \: \div F , \quad x \in {\cal O},\\
 \div u  & \: =  
\: g, \quad x \in {\cal O}, \\
 u\vert_{\pa {\cal O}}  & \: =  \: 0,
\end{aligned}
\right.
\end{equation}
when ${\cal O}$ is a bounded domain with $C^{1,1}$ boundary.
{\em Such inequality is no longer valid when $\pa {\cal O}$ has weaker
  regularity}. Only
 the interior estimate 
\begin{equation} \label{interior}
\| \na u \|_{H^1(K)} + \| p \|_{H^1(K)/\R} \: \le C(K) \, \left( \| F
\|_{H^1({\cal O})} \: + \:   \| g \|_{H^1({\cal O})} \right) 
\end{equation}
is satisfied, where $K$ is any relatively compact open subset of ${\cal
  O}$. As regards well-posedness issues, this interior bound is not
sufficient. We will need a control up to the boundary, given by the
following:

\begin{proposition} \label{StokesBMO} 
Assume that ${\cal O}$ has a $C^{1,\alpha}$ boundary, $\: 0 < \alpha \le
1$. Assume also that  
$$ F \in  L^2({\cal O})  \cap \bmo({\cal O}), \quad g \in L^2({\cal O})
\cap \bmo({\cal O}). $$ 
Then, the weak solution $(u,p)$ of \eqref{Stokes} satisfies 
\begin{equation} \label{BMO} 
\| \left( \na u, p \right)  \|_{\bmo(\O)}  \: \le C \, \left(  \| \left( F, g \right) 
\|_{\bmo(\O)} \: + \: \| \left( F, g \right) \|_{L^2(\cal O)} \right).
\end{equation}
\end{proposition}
We remind that $\bmo({\cal O})$ is the set of functions $f \in L^1({\cal O})$ such that 
$$ \sup_{B} \frac{1}{|B|} \, \int_{B} | f(x) -  \overline{f}_{B} |  \, dx
 \: < \: +\infty,  
 \quad  \overline{f}_{B} \: = \:  \frac{1}{|B|} \int_B f(x) dx,   $$
 where the supremum is taken over all the open balls $B$ of $\O$, that is
 all the intersections of $\O$ with open disks.  Note that the application  
$$ || f ||_{\bmo(\O)} \: := \:  \sup_{B} \frac{1}{|B|} \, \int_{B} | f(x) -
 \overline{f}_{B} |  \, dx   $$ 
 defines only a semi-norm, as  it is invariant by the addition of
 constants. An easy remark is that   $\bmo(\O)$ is also characterized by   
$$ \sup_{B}  \, \inf_m \left( \frac{1}{|B|} \, \int_{B} | f(x) -  m |  \,
 dx  \: \right) < \: +\infty, $$   
where the infimum is taken over all real constants, providing an
 equivalent semi-norm.
 Note that $f \in \bmo({\cal O})$ if and only if 
 $\tilde{f} \in \bmo(\R^2)$, where $\tilde{f}$ is the extension of $f$ by
 zero. Hence, standard results for the whole space  apply  directly to our
 setting.  
 For instance,  $ f \in \bmo(\O)$ belongs to $L^p(\O)$ for any finite $p$, and 
  $$ \sup_{B} \left( \frac{1}{|B|} \, \int_{B} | f(x) -  \overline{f}_{B}
 |^p  \, dx  \: \right)^{1/p} < \: +\infty, $$   
 this expression defining again a semi-norm which is equivalent to the previous
 one.  We also remind the Sobolev imbedding in dimension 2: 
\begin{equation} \label{imbedding} 
 H^1(\O) \hookrightarrow \bmo(\O), \quad \| f \|_{\bmo(\O)} \: \le \: C \,
 \| f \|_{H^1(\O)}.  
 \end{equation}
which is simply deduced from Poincar\'e inequality. Finally, we remind the
interpolation inequality: for all $\theta \in (0,1)$, for all $1 \le p, q <
+\infty$ with $(1-\theta) q = p$  
\begin{equation} \label{interpolation} 
  \| f \|_{L^q(\O)} \: \le \: C \, \| f \|_{L^p(\O)}^{1-\theta} \,  \| f
  \|_{\bmo(\O)}^\theta  
\end{equation}
 We refer to \cite{Fefferman&Stein72,Journe83} for exhaustive study of the space $\bmo$. 

\medskip
{\bf Proof of the proposition.} 
In the case of the whole space $\O = \R^2$, the estimate  \eqref{BMO}
follows from the continuity of the Riesz transform on $\bmo$. In the case
of a $C^{1,\alpha}$ bounded domain,  it is connected  to  H\"older theory
for elliptic systems. Such theory  has been of course widely considered,
from various perspectives: see
\cite{ADN59,Benmomo&Robbiano91,Giaquinta83,Giaquinta&Modica82}  
for some examples. Although a $\bmo$ estimate like \eqref{BMO} may be 
part of the folklore of this domain, we could not find a proper
reference for it. For the sake of completeness, we  give here  the main
steps of (one possible) proof. The last step of the proof relies on ideas
of Giaquinta et Modica, used to establish  
H\"older estimates for the Stokes system with Neumann boundary condition
\cite{Giaquinta&Modica82}.   

\medskip
We start with a simple remark, to be used implicitly throughout the sequel:
any $f \in L^2(U)$, $U$ open set, can be written $f = \div F$,  where $F
\in H^1(U)$ satisfies  
\begin{equation} \label{controlBMO} 
\| F \|_{\bmo(U)} \: \le \: C  \| F \|_{H^1(U)} \:  \le \: C' \, \| f \|_{L^2(U)} 
\end{equation}
This will allow  to keep the source term in divergence form as we apply
transformations to the Stokes system.

\medskip
Let $(u, p)$ be the weak solution of \eqref{Stokes}, where $p$ is normalized such that 
$ \int_{\cal O} p =  0$. Standard energy estimates yield 
\begin{equation} \label{standardestimate} 
 \| u \|_{H^1(\O)}  \: + \: \| p \|_{L^2(\O)} \: \le \: C \, \left( \| F
 \|_{L^2(\O)} \: + \:  
\| g \|_{L^2(\O)} \right).
\end{equation}

\medskip
{\em Step 1 : Localization.} Let  $\tilde{\O}^{i}  \Subset \O^i$, $\,
i=1\dots N$ a covering  of $\overline{\O}$ by open sets.  Let $\psi^i$ a smooth
function with compact support in $\O^i$, such that $\, \psi^i = 1$ on
$\tilde{\O}^{i}$.  The functions  
$$u^i \: :=  \: \psi^i \, u,  \quad p^i \: :=  \: \psi^i \, p$$ 
satisfy  
\begin{equation} \label{Stokesi}
\left\{
\begin{aligned}
 -\Delta u^i  + \na p^i   \: =  &  \: \div\left(\psi^i F\right) \: - \: F \, \na \psi^i
\: - \; 2 (\na u)^t \na \psi^i
  - \Delta \psi^i \, u  \: + \: p \na \psi^i \: \\  
  \: :=  & \: \div F^i, \quad x \in {\cal O}^i,\\
 \div u   \: =  & \: g \, \psi^i  \; + \:  \na \psi^i \cdot u 
   := \:  g^i ,  \quad x \in {\cal O}^i, \\
 u\vert_{\pa {\cal O}^i}  \: = &  \: 0.
\end{aligned}
\right.
\end{equation}
By \eqref{standardestimate}, the $L^2 \cap \bmo$ norms
of $F^i$ and $g^i$ are controlled by the $L^2 \cap \bmo$ norms  of $F$ and
$g$.  Thus, we can restrict ourselves to  a
subdomain, that is establish \eqref{BMO} with $\O^i$ instead of  $\O$.

\medskip
{\em Step 2 : Local coordinates.} If $\O^i$ does not intersect the boundary
of $\O$,  the estimate follows from the interior regularity
\eqref{interior}. If $\O^i$ intersects the boundary,  we can assume with no
loss of generality that  it  is a local chart: there exists a $C^{1,
  \alpha}$ diffeomorphism  
$$\chi : \O^i \mapsto D(0,R), \quad \chi\left( \O^i  \cap \pa \O \right)
\: = \: (-R,R) \times  
\{ 0 \}, \quad \chi\left( \O^i  \cap \O \right)   \: = \: D^+(0,R), $$
where $D^+(0,R)$ is the upper half disk of radius $R$ centered at the
origin. We define new fields $v$, $q$, $F'$, $g'$ by the relations 
$$ u^i(x) \: :=  \: v(\chi(x)),  \quad p^i(x) \: := \: q(\chi(x)), \quad
F^i(x) \: = \: F'(\chi(x)), \quad g^i(x) \: =  \: \dfrac{g'(\chi(x))}{\text{det}(\nabla \chi)}. $$  
They satisfy 
\begin{equation*} 
\left\{
\begin{aligned}
 -\div (A \na  v) + \div(B q)   & \: =  \: \div(B F'), \quad x \in D^+(0,R),\\
 B :  \na v  & \: = \; g' , \hspace{1.6cm} x \in D^+(0,R), \\
 v\vert_{\pa D^+(0,R)}  & \: =  \: 0.
\end{aligned}
\right.
\end{equation*}
where 
$$ A =  \frac{1}{\det(\na \chi)} \, (\na \chi)^t \, \na \chi, \quad B =
\frac{1}{\det(\na \chi)} (\na \chi)^t. $$ 
Note that $A$ is uniformly elliptic over $D(0,R)$, and that $A$, $B$ have
$C^{0,\alpha}$  coefficients. As usual, $\: (\div M)_i := \pa_j M_{ji}$, and
$\: M : N = M_{ij} \, N_{ij}$ for any 2x2 matrices $M, N$. 
 
\medskip
{\em Step 3 : Frozen coefficients}. 
We write the previous system as 
\begin{equation} \label{stokesfrozen} 
\left\{
\begin{aligned}
 -\div (A(0) \na  v) + \div(B(0) q)   & \: =  \: \div(\tilde{F}), \quad x \in D^+(0,R),\\
 B(0) : \na v  & \: = \; \tilde{g}  
, \hspace{1.3cm} x \in D^+(0,R), \\
 v\vert_{\pa D^+(0,R)}  & \: =  \: 0.
\end{aligned}
\right.
\end{equation}
where 
$$ \tilde{F} \: : = \: B F' \: - \: (A(0) - A(x)) \na v + (B(0)-B(x)) q,
\quad \tilde{g} \: := \: g' + (B(0)-B(x)) : \na v. $$   
Let us  assume for a while that $\tilde{F} \in L^2 \cap \bmo$, $\tilde{g}
\in L^2 \cap \bmo$, and that the estimate   
\begin{equation} \label{BMObis}  
\| \left(\na v, q\right)  \|_{\bmo(D^+(0,R))}   \: \le {\cal C} \,   \Bigl(
\| ( \tilde{F}, \tilde{g} ) 
\|_{\bmo(D^+(0,R))} \: + \:   \|( \tilde{F}, \tilde{g} )  \|_{L^2(D^+(0,R))} \Bigr).
\end{equation}
holds. A simple scaling argument shows that the constant ${\cal C}$ can  be
chosen independently of the radius $R$. We now state the following  {\it a
  priori} estimate: there exists a universal constant ${\cal C'}$, and
$\eps(R)$ going to zero with $R$ such that  
\begin{equation} \label{estimateFg}
\begin{aligned}
&  \| (\tilde{F}, \tilde{g}) \|_{\bmo(D^+(0,R))} \: +  \: \| ( \tilde{F},
 \tilde{g}) \|_{L^2(D^+(0,R))} \:  \le \: \eps(R) \, \| (\na v,q)
 \|_{\bmo(D^+(0,R))}  \\  
 & + \:  {\cal C}' \left( \| (F', g' ) \|_{\bmo(D^+(0,R))} \: +  \: \| ( F', g' ) 
 \|_{L^2(D^+(0,R))} \,  + \,  \| (\na v,q) \|_{L^2(D^+(0,R))} \right)  
\end{aligned}
\end{equation}
For the sake of brevity, we focus on the $\bmo$ bound, as the $L^2$ bound
is straightforward. More precisely, we just show how to bound  
$\displaystyle  \|  (A(0) - A(x)) \na v \|_{\bmo}, $
because the other terms composing $\tilde{F}$ and $\tilde{g}$ can be
handled along the same lines. As emphasized at the beginning of the
section, we need to control 
$$ I_B \: := \: \frac{1}{|B|} \int_B \bigl| \left(A(x)-A(0)\right) \, \na v(x)
\, -  \,  c \bigr| \, dx $$ 
for any ball $B$ of $D^+(0,R)$ and some constant vector $c$ (possibly depending on
$B$). Let $r$ be the diameter of $B$ and $x_0$ a point in $B$. We choose
$\displaystyle c
= \left( A(x_0) - A(0) \right)  \overline{(\na v)}_B$. We get  
\begin{equation} \label{calculIB}
\begin{aligned}
I_B \: &  \le \:  \frac{C}{r^2}  \, \left( \int_B |A(x) - A(x_0)| \,
\overline{(\na v)}_B \, dx \: + \:  
\int_B \,  |A(x) - A(0)|  \, | \na v(x) - \overline{(\na v)}_B  | \, dx \right)  \\
& \le \:   C' \left(  r^{\alpha-2}  \, \int_B |\na v(x)| dx \: + \: R^\alpha \,
\| \na v \|_{\bmo(D^+(0,R))} \right) \\ 
& \le \:   C' \left( r^{\alpha-2+2/q} \, \|\na  v \|_{L^p(D^+(0,R))} \: + \:
R^\alpha \, \|\na  v \|_{\bmo(D^+(0,R))} \right) 
\end{aligned}
\end{equation}
for any finite conjugate exponents $p,q$, {\it i.e.} $p^{-1} + q^{-1} =
1$. We choose $q$ close enough to 1 so that  $\alpha-2+2/q  > 0$. Together
with the interpolation inequality \eqref{interpolation}, we deduce that  
\begin{align*}
 I_B \: & \le \: C'' \, \left( R^{\alpha-2 + 2q} \| \na v
 \|_{L^2(D^+(0,R))}^{2/p} \, \|  \na  v  
 \|_{\bmo(D^+(0,R))}^{1-2/p} \: + \:  R^\alpha \, \| \na v \|_{\bmo(D^+(0,R))} \right) \\
 & \le \:  C_1 R^{\gamma} \, \|\na v \|_{\bmo(D^+(0,R))}  + C_2   
\| \na v \|_{L^2(D^+(0,R))}
 \end{align*}
 for some universal positive constants $\gamma$, $C_1, C_2$. The estimate
 \eqref{estimateFg} follows.  
 
 \medskip
Note that estimates \eqref{BMObis} and  \eqref{estimateFg} yield the bound
\eqref{BMO}.  Indeed,    
up to take smaller $R$, that is up to refine the  covering  of open sets $\O^i$,
we can  assume that $\eps(R) \le 1/(2 {\cal C})$. Hence, combining
\eqref{estimateFg}-\eqref{BMObis}, we obtain  
\begin{align*}
  \frac{1}{2} \| \left(\na v, q\right)  \|_{\bmo(D^+(0,R))}   & \:  \le \:
  \frac{1}{2} \,   \| (\na v,q) \|_{L^2(D^+(0,R))} \\
&  + \:  {\cal C} {\cal C'}\,   \Bigl(  \| ( F', g')
\|_{\bmo(D^+(0,R))} \: + \:   \|( F', g') )  \|_{L^2(D^+(0,R))} \Bigr). 
\end{align*}
Then, it is well-known that  $L^2$, $H^1$ and $\bmo$ norms are preserved by
$C^1$ diffeomorphisms. This allows to bound the right-hand side of the
previous inequality:
\begin{align*}
 & \frac{1}{2} \,   \| (\na v,q) \|_{L^2(D^+(0,R))} 
  + \:  {\cal C} {\cal C'}\,   \Bigl(  \| ( F', g')
\|_{\bmo(D^+(0,R))} \: + \:   \|( F', g') )  \|_{L^2(D^+(0,R))} \Bigr) \\
& \le \: C \, \left( \| u^i \|_{H^1({\cal O}^i \cap {\cal O})}  \: + \: 
\| p^i \|_{L^2({\cal O}^i \cap {\cal O})} \: 
+ \: \| (F^i,g^i) \|_{L^2({\cal O}^i \cap {\cal O})} \: + \: 
\| (F^i,g^i) \|_{\bmo({\cal O}^i \cap {\cal O})} \right)\\
& \le \: C'  \left( \| (F,g) \|_{L^2({\cal O})} \: + \: 
\| (F,g) \|_{\bmo({\cal O})} \right) 
\end{align*}
where the last line involves  the basic estimate
\eqref{standardestimate}. As regards the left-hand side, we obtain the
lower bound
$$ \| p^i \|_{\bmo({\cal O}^i \cap {\cal O})} \: \le \: C \: \| q \|_{\bmo(D^+(0,R))}   $$
and along the  lines of \eqref{calculIB}
$$
\begin{array}{rcl}
 \| \na u^i \|_{\bmo({\cal O}^i \cap {\cal O})} &\: = \:& \| \na \chi \na
v(\chi(\cdot)) \|_{\bmo({\cal O}^i \cap {\cal O})}, \vspace{.1cm}\\
& \: \le \: & C \left( \|
\na v \|_{\bmo(D^+(0,R))} \: + \: \| \na v \|_{L^2(D^+(0,R))}\right).
\end{array}
$$
This altogether implies \eqref{BMO}.

\medskip
We stress that \eqref{BMObis} and \eqref{estimateFg} are only {\it a
  priori} estimates:  $\na u,p$, and therefore $\na v,q$ are only supposed
to be in $L^2$, and not in $\bmo$. Nevertheless,  regularizing the
coefficients of $A$ and $B$, establishing the same estimates for the
regularized problem and passing to the limit  allows to show that the weak
solutions are indeed in $\bmo$ and that the inequality holds. As this
regularization argument is very classical, we leave it to the reader.  

\medskip
{\em Step 4: $\bmo$ estimate for the Stokes system in a half-disk}.  
The final step of the proof is to derive the estimate \eqref{BMObis}  for
the system \eqref{stokesfrozen}.  By the {\em reverse} change of variables: 
$$ x \mapsto  \: \left( \na \chi(0)^t  \right)^{-1} \, x $$
we can assume that $A(0) = B(0) = I_2$ is the identity matrix. By this
linear mapping, the domain $D^+(0,R)$ turns into the  intersection of a
half-plane and an ellipse, say $E^+$. As all the vector fields involved are
compactly supported in $E^+$,  this Stokes system with Dirichlet boundary
conditions still holds in any half-disk containing $E^+$. As this system
is rotationally invariant, we can furthermore assume the half-disk to be
$D^+(0,R')$ for some large enough $R'$. Finally, as the estimate
\eqref{BMObis} is  invariant by the dilations  $x \mapsto R'x$, we can
consider the case $R'=1$. Eventually, we only have to establish the
inequality  
$$ \| (\na u,p)  \|_{\bmo(D^+(0,1))} \: \le  \: C \, \left(   \| (F,g)
\|_{\bmo(D^+(0,1))} \: + \:  \| (F,g) \|_{L^2(D^+(0,1))}  \right) $$ 
for the system 
\begin{equation*} 
\left\{
\begin{aligned}
 -\Delta u + \na p   & \: =  \: \div F, \quad x \in D^+(0,1),\\
\div u  & \: = \; g  
, \hspace{1cm} x \in D^+(0,1), \\
 u\vert_{\pa D^+(0,1)}  & \: =  \: 0.
\end{aligned}
\right.
\end{equation*}
We remind that, if $p$ is chosen such that $\displaystyle \int_{D^+(0,1)}
\!\!\!\!\! \! p = 0$, we already have the $L^2$ estimate 
$$ \| (\na u,p)  \|_{L^2(D^+(0,1))} \: \le  \:  C \:  \| (F,g) \|_{L^2(D^+(0,1))}.  $$

We shall rely on  ideas  of Giaquinta and Modica, who prove in article
\cite{Giaquinta&Modica82} a H\"older  estimate for the Stokes equation with Neumann type
boundary conditions. Let $0 < \rho \le R \le 2$, and $y$ in $D^+(0,1)$.  We
will denote  
$$B(y, \rho) \: := \: D(y,\rho) \cap D^+(0,1), \quad  \overline{f}_{y,\rho}
\: := \:  \frac{1}{\left|B(y,\rho)\right|}  \int_{B(y,\rho)} f, \quad
\forall \, 0 < \rho \le R.$$    
We decompose  $u=v+w$, $p=q+r$, where  $(v,q)$ solves  
\begin{equation*} 
\left\{
\begin{aligned}
 -\Delta v + \na q   & \: =  \: \div F, \hspace{0.8cm}  x \in B(y,R),\\
\div  v  & \: = \; g  - \overline{g}_{y,R}, \quad 
 x \in B(y,R), \\
 v\vert_{\pa B(y,R)}  & \: =  \: 0.
\end{aligned}
\right.
\end{equation*}
and $(w,r)$ solves 
\begin{equation*} 
\left\{
\begin{aligned}
 -\Delta w + \na r   & \: =  \: 0, \hspace{0.8cm} x \in B(y,R),\\
\div w  & \: = \;    \overline{g}_{y,R}, \quad x \in B(y,R), \\
 w\vert_{\pa B(y,R)}  & \: =  \: u\vert_{\pa B(y,R)}.
\end{aligned}
\right.
\end{equation*}

We must first derive an estimate on $v$ and $q$. We state without proof the
well-known inequality  (see \cite{Sohr01}) 
\begin{equation} \label{estimq} 
\| q - \overline{q}_{y,\rho} \|_{L^2(B(y,\rho))} \: \le \: C \| \na v
\|_{L^2(B(y,\rho))},  
\end{equation} 
where $C$ does not depend on $\rho$ by a simple scaling argument. 
Then, a standard energy estimate yields 
\begin{equation*}
 \int_{B(y,R)} |\na v |^2 \:   = \: -\int_{B(y,R)}  \left( F -
 \overline{F}_{y,R} \right)  \cdot \na v  \: + \: \int_{B(y,R)} \left( g -
 \overline{g}_{y,R} \right)  \cdot \left( q -
 \overline{q}_{y,R} \right) 
\end{equation*} 
which combined with \eqref{estimq} yields
\begin{equation} \label{estimv}
 \| \na v\|_{L^2(B(y,R))} \: \le \: C \left(  \| F -  \overline{F}_{y,R}
 \|_{L^2(B(y,R))} \: + \:  
 \| g -  \overline{g}_{y,R} \|_{L^2(B(y,R)} \right). 
\end{equation}

We now wish to obtain an estimate on $w$ and $r$. As for $q$, the pressure $r$ satisfies
\begin{equation} \label{estimr} 
\| r - \overline{r}_{y,\rho} \|_{L^2(B(y,\rho))} \: \le \: C \| \na w
\|_{L^2(B(y,\rho))},  
\end{equation}

As regards $w$, we want to show the estimate 
\begin{equation} \label{estimw}
 \| \na w - \overline{(\na w)}_{y,\rho} \|_{L^2(B(y,\rho))} \: \le \: C
 \frac{\rho^2}{R^2} \| \na w - \overline{(\na w)}_{y,R} \|_{L^2(B(y,R))} 
 \end{equation}
where  $C$ does not depend on $\rho$ or $R$. At first, up to replace $w$ by
$w - x_2 \left( \begin{smallmatrix} \overline{(\pa_2 w_1)}_{y,R} \\
  \overline{g}_{y,R}\end{smallmatrix} \right)$, which would still be zero at the
flat part of the  boundary $\pa B(y,R) \cap \{x_2=0\}$, and would still
satisfy \eqref{estimw}, we can assume that  
 \begin{equation*}
\overline{(\pa_2 w_1)}_{y,R}  \: = \: \overline{g}_{y,R} \: = \: 0.
 \end{equation*}
If $\rho > R/2$, inequality \eqref{estimw} is trivially satisfied. If $\rho
< R/2$, there are two  
cases. 

\medskip
{\em If $B(y,R) \subset \{ x_2 > 0\} $}, the ball $B(y,R)$ does not
intersect the boundary of $D^+(0,1)$. We can use the interior estimate
provided by Giaquinta and Modica in \cite{Giaquinta&Modica82}: we can apply
proposition 1.9, 
estimate (1.14) to the derivatives of $w$, which are still solutions of the
Stokes equation, and this yields  exactly \eqref{estimw}. 

\medskip
{\em If  $B(y,R) \cap \{ x_2 =  0\} \neq \emptyset$}, we write 
\begin{align*}
\| \na w - \overline{(\na w)}_{y,\rho} \|_{L^2(B(y,\rho))} \: & \le \: C \,
\rho \, \| \na^2 w \|_{L^2(B(y,\rho))} \\  
& \le \: C \, \rho^2 \, \| \na^2 w \|_{L^\infty(B(y,\rho))}  \: \le \: C
\rho^2 \, \| \na^2 w \|_{L^\infty(B(y,R/2))} \\ 
& \le \: C(R) \, \rho^2 \| \na w \|_{L^2(B(y,R))} 
\end{align*}
Note that the first inequality is simply Poincar\'e's inequality, whereas
the last one stems from classical regularity results for the
Stokes operator. Simple scaling considerations give the bound $C(R) \le
C/R^2$ for some constant $C$ that does not depend on $R$. To prove
\eqref{estimw}, it is therefore enough to show that: for any solution $w$
of the Stokes equation  
$$
\left\{
\begin{aligned}
 -\Delta w + \na p   & \: =  \: 0, \hspace{1cm} x \in  B(y,R),\\
\div w  & \: = \; 0  
, \hspace{1cm} x \in B(y,R), \\
\end{aligned}
\right.
$$
satisfying moreover
\begin{equation}  \label{assumew} 
\overline{(\pa_2 w_1)}_{y,R}  \: =  \: 0, \quad  w = 0 \: \mbox{ on } \quad
\pa B(y,R) \cap \{ x_2 =  0\}  
\end{equation}
we have 
\begin{equation} \label{estimw2}
 \| \na w \|_{L^2(B(y,R))} \: \le \: C \, \| \na w - \overline{(\na w)}_{y,R}
 \|_{L^2(B(y,R))} \quad \mbox{ if } \pa B(y,R) \cap \{ x_2 =  0\} \neq
 \emptyset.  
 \end{equation} 
 Again,  the constant $C$ in the r.h.s can  be chosen  independently  of $R$. 
 
 \medskip
 If inequality \eqref{estimw2} were not to be satisfied, one could find
 solutions $w^n$ satisfying \eqref{assumew}, and such that  
 $$ \|   \na w^n  \|_{L^2(B(y,R))} = 1, \quad  \| \na w -
   \overline{(\na w^n)}_{y,R} \|_{L^2(B(y,R))} \xrightarrow[n \rightarrow
   +\infty]{} 0.  $$ 
   From the first equality,  up to a subsequence,  $w^n
   \rightarrow w$ weakly in $H^1(B(y,R))$. This implies the convergence of
   the averages  
$\overline{(\na w^n)}_{y,R} \rightarrow \overline{(\na w)}_{y ,R}$. Moreover, by
   standard ellipticity properties of the Stokes operator, we have  
$$ \| \na^2 w^n \|_{L^2(B)} \: \le \: C(R,B), \quad \forall \,B \Subset B(y,R) $$
so that $w^n \rightarrow w$ strongly in  $H^1(B)$.
 Hence, we obtain, as $n \rightarrow +\infty$, 
 $$ \na w \: = \:   \overline{(\na w)}_{y,R} \quad \mbox{ on }  B. $$

 From the second condition in \eqref{assumew} and the divergence-free condition, we get
$$ \overline{(\pa_1 w)}_{y,R} = 0, \quad \overline{(\pa_2 w_2)}_{y,R} = 0. $$
Moreover, by the first condition in \eqref{assumew}, we also have
$\overline{(\pa_2 w_1)}_{y,R}  =  0$. 
Hence, $\overline{(\na w)}_{y,R} = 0$, and $\na w=0$ in any subset $B$
relatively compact in $B(y,R)$. Thus,  $\na w = 0$ on $B(y,R)$ which
contradicts the assumption that its $L^2$ norm is 1.  

\medskip
 This last argument leads to the desired inequality \eqref{estimw2} on
 \eqref{estimw}. Combining  \eqref{estimv} and \eqref{estimw}, we obtain  
 \begin{align*}
 &  \| \na u \: - \:   \overline{(\na u)}_{y,\rho} \|_{L^2(B(y,\rho))}\:
  \le \:  \| \na w \: - \:   \overline{(\na w)}_{y,\rho} \|_{L^2(B(y,\rho))}
\: + \:  \| \na v \: - \:   \overline{(\na v)}_{y,\rho} \|_{L^2(B(y,\rho))}
\\
&   \le \: C \, \left( \frac{\rho^2}{R^2} \,  \| \na w \: - \:   \overline{(\na
    w)}_{y,R} \|_{L^2(B(y,R))} \: + \: \| \na v \|_{L^2(B(y,R))} \right) \\
&  \le \: C' \, \biggl( \frac{\rho^2}{R^2} \, 
  \| \na u \: - \:   \overline{(\na u)}_{y,R} \|_{L^2(B(y,R))} \: + \:    \|
   F -  \overline{F}_{y,R} \|_{L^2(B(y,R))} \: + \:  \| g -
   \overline{g}_{y,R} \|_{L^2(B(y,R))} \biggr)  \\ 
 &  
 \le \:  C'' \, \left( \frac{\rho^2}{R^2} \, 
  \| \na u \: - \:   \overline{(\na u)}_{y,R} \|_{L^2(B(y,R))} \: + \:   \|
   (F,g) \|_{\bmo(D^+(0,1))}  
  \,  R^2  \right)
\end{align*}
 We use lemma 0.6 of \cite{Giaquinta&Modica82} to conclude that 
 $$ \| \na u \: - \:   \overline{(\na u)}_{y,\rho} \|_{L^2(B(y,\rho))} \: \le
 \: C \, \| (F,g) \|_{\bmo(D^+(0,1))} \,  \rho^2, $$ 
which provides the BMO control of $\na u$. The BMO control of the pressure
$p$ then follows from \eqref{estimq}, \eqref{estimr}. This ends the proof.   

\bigskip
In the next section, we will use this proposition  to show well-posedness
of  the PDE's system \eqref{NS}-\eqref{contin2}. Before that, we state a
regularity result of Sobolev type  for the Stokes system in $C^{1,\alpha}$
domains. It will allow to give a meaning  in the trace sense to the stress
tensor at the solid  boundary $(\pa_n u - p \, n)\vert_{\pa S^i}$.  
\begin{proposition} \label{StokesSob} 
Assume that ${\cal O}$ has a $C^{1,\alpha}$ boundary, $\: 0 < \alpha \le
1$. Let $s,\tau$ such that $s < \alpha$ and $s \le 2/\tau$.  Assume that  
$$ F \in  L^2({\cal O})  \cap W^{s,\tau}({\cal O}),  \quad g \in L^2({\cal O})
\cap W^{s,\tau}({\cal O}). $$ 
Then, the weak solution $(u,p)$ of \eqref{Stokes} satisfies 
\begin{equation} \label{Wsp} 
\| \left( \na u, p \right)  \|_{W^{s,\tau}(\O)}  \: \le C \, \left(  \|
  \left( F, g \right) v 
\|_{W^{s,\tau}(\O)} \: + \: \| \left( F, g \right) \|_{L^2(O)} \right).
\end{equation}
\end{proposition}
We remind that for all $0 < s < 1$, the fractional Sobolev space
$W^{s,\tau}(\O)$ is the set of measurable functions $u$ such that  
$$ \| u \|_{W^{s,\tau}(\O)} \: := \:  \left( \int\int_{\O \times \O}
 \frac{|u(x) - u(y)|^\tau}{|x-y|^{2+s\tau}} \, dxdy \right)^{1/\tau} \: < \: +
 \infty. $$ 
 and this last expression makes it a Banach space. The assumption  $s \le
 2/\tau$ in the proposition ensures the continuous imbedding 
 \begin{equation} \label{imbedding2}
 H^1(\O) \hookrightarrow W^{s,\tau}(\O), \quad \| f \|_{W^{s,\tau}(\O)} \: \le \:
 C \, \| f \|_{H^1(\O)}.  
 \end{equation}
Similarly, the constraint $s < \alpha$ is such that $C^{0,\alpha}(\O)
\hookrightarrow   W^{s,\tau}(\O)$.  

\medskip
{\bf Sketch of proof of the proposition.} The proof of the Sobolev estimate
\eqref{Wsp} mimics the proof of the $\bmo$ estimate \eqref{BMO}, so that we
only quote the few changes to be made.  
 
\medskip
Steps 1 and 2 (localization and use of local coordinates) remain the same,
up to the replacement of  $\bmo$ by $W^{s,\tau}$ in every argument.    

\medskip
In step 3, the only change is in the derivation of 
\begin{equation*} 
\begin{aligned}
&  \| (\tilde{F}, \tilde{g}) \|_{W^{s,\tau}(D^+(0,R))} \: +  \: \| (
 \tilde{F}, \tilde{g}) \|_{L^2(D^+(0,R))} \:  \le \: \eps(R) \, \| (\na v,q)
 \|_{W^{s,\tau}(D^+(0,R))}  \\  
 & + \:  {\cal C}' \left( \| (F', g' ) \|_{W^{s,\tau}(D^+(0,R))} \: +  \: \| ( F', g' ) 
 \|_{L^2(D^+(0,R))} \,  + \,  \| (\na v,q) \|_{L^2(D^+(0,R))} \right) 
\end{aligned} 
\end{equation*}
which substitutes to \eqref{estimateFg}. Again,  we just show how to bound
$\| (A(0) - A(x) ) \na v \|_{W^{s,\tau}}$, as all other terms that compose
$\tilde{F}$ and $g$ are treated in the same manner. We write 
\begin{align*}
& \int \int_{D^+(0,R) \times D^+(0,R)} \frac{\left| (A(0) - A(x) ) \na v(x)
    - (A(0) - A(y)) \na v(y)\right|^\tau}{|x-y|^{2+s\tau}} \, dxdy \\ 
&  \le \: C \bigl(  \int \int_{D^+(0,R) \times D^+(0,R)}
  |A(0) - A(x)|^\tau \frac{|\na v(x) - \na v(y)|^\tau}{|x-y|^{2+s\tau}} \, dxdy \:  \\
&  \quad + \:    \int \int_{D^+(0,R) \times D^+(0,R)}  |\na v(y)|^\tau
  \frac{|A(x)-A(y)|^\tau}{|x-y|^{2+s\tau}} \, dxdy\bigr)\\ 
  &  \le \: C \: \left(   \, \| \na v \|_{W^{s,\tau}(D^+(0,R))}^\tau  \, \sup_{x\in
    D^+(0,R)} |A(x)-A(0)|^\tau \: + \:    R^{\tau(\alpha-s)}  \, 
   \| \na v \|_{L^\tau(D^+(0,R))}^\tau \, \right) \\
  & \le \:  \eps'(R) \,   \|\na v \|_{W^{s,\tau}(D^+(0,R))}^\tau. 
  \end{align*}
 which allows to conclude as in the previous proof. 
 
 \medskip
Step 4, that is  the $W^{s,\tau}$ estimate, $0 < s < 1$, for the  Stokes
equation in a half-disk, follows from a simple interpolation of similar
inequalities for $W^{0,\tau}$ and $W^{1,\tau}.$

\section{Strong solutions}
This section is devoted to the proof of theorem \ref{theo1}. {\em Broadly, we shall prove
existence and uniqueness of strong solutions as long as the distance
between solid boundaries $\delta(t)$ satisfies $\delta(t) > \delta_0$,      
where $\delta_0 > 0$ is arbitrary}. The fact that strong solutions can not
exist after collision will be discussed eventually.  This altogether will
of course imply the result. In what follows, constants will depend
implicitly on $\delta_0$.  

\medskip
We treat separately the existence and uniqueness parts. The existence result
follows the lines of \cite{Desjardins&Esteban99}, whereas the uniqueness result is
inspired by \cite{Takahashi03bis}. We thus rely substantially on these articles, and
put the stress only on the changes due to our not so regular $C^{1,\alpha}$
boundaries.  

\medskip
Our (refined) existence result reads:
\begin{proposition} {\bf (Existence of strong solutions)} 

\smallskip 
Let $\delta_0 > 0$, $\displaystyle 
\: v_0 \in H^1_0(\Omega), \quad \rho^i_0 \, D(v_0) = 0, \: \forall \, i, 
 \quad  \: f \in L^2((0,T) \times \Omega),  \:  \forall \, T >0.$ 
 Assume \eqref{nocontact}, and
$$\displaystyle \pa \Omega \in C^{1, \alpha}, \quad \pa S^i_0 \in C^{1,
   \alpha}, 
\quad \forall \, i, \quad 0 < \alpha \le 1.$$
 Then there exists a  strong  solution on $(0,T)$ for some $T >
   0$. Moreover, one of the following alternatives holds true: 
 \begin{description}
\item[i)]  One can take $T$ arbitrarily large and $\delta(t) > \delta_0$ for
  all $t \le T$.  
\item[ii)]  One can take $T$ such that $\delta(t) > \delta_0$ for all $t < T$
  and $ \:\lim_{t \rightarrow T} \delta(t) = \delta_0$.  
\end{description} 
In both cases, the strong solution has the additional regularity
 $$ \int_0^T \| \na v(t) \|_{\bmo(F(t))}^2   dt  \: + \:   \int_0^T \| q
 \|_{\bmo(F(t))}^2   dt  < +\infty $$ 
 and  
 $$ \int_0^T \| \na v(t) \|_{H^1(F_{\eps}(t))}^2   dt  \: + \:   \int_0^T \| q
 \|_{H^1(F_{\eps}(t))/\R}^2   dt  < +\infty, $$ 
 where $q$ is the corresponding pressure field, and 
$$
F_{\eps}(t) :=\{x \in F(t) \; \mathrm{s.t.} \; \mathrm{dist}(x,\partial F(t)) >
\eps\}, \quad \eps >0.
$$  
  \end{proposition}

\medskip
{\bf Proof of the proposition}.  Following \cite{Desjardins&Esteban99}, we
establish {\it a priori} estimates for a sufficiently smooth solution
$(v,q)$ on $(0,T)$,  s.t.  $\delta(t)>\delta_0$ for all $t<T$.

\medskip
 We first take $\varphi = v$ as a test function, which yields the standard
 energy inequality  
\begin{equation} \label{apriori}
\| v \|_{L^\infty(0,T; \, L^2(\Omega))} + \| v \|_{L^2(0,T; H^1_0(\Omega))}
\: \le \: C(T) \left( \| v_0 \|_{L^2(\Omega)} \: + \: \| f \|_{L^2((0,T)
    \times \Omega)} \right) 
\end{equation}
where $C(T)$ is an increasing function of $T$.  Then, we take
$\varphi=\pa_t  v$ as a test function, which yields 
\begin{equation} \label{apriori1}
\int_0^t \int_\Omega |\pa_t v |^2 \: + \: \mu \, \int_\Omega | D(v)(t) |^2 \: \le \: 
C \left( \int_\Omega | D(v_0) |^2  \: + \: \int_0^t \int_\Omega  | f |^2 \: + \: 
 \int_0^t \int_\Omega  | v \cdot \na v  |^2 \right)
\end{equation}
Note that the l.h.s. in \eqref{apriori}, {\em resp.} \eqref{apriori1}
controls the $L^\infty \cap L^2$ norm  of $v^i, \omega^i$, {\em resp.} the
$L^2$ norm of $\dot{v}^i, \dot{\omega}^i$. We now use Stokes regularity to
bound the last term in \eqref{apriori1}. 

\medskip
The Navier-Stokes equation for the fluid part can be written
\begin{equation} \label{NSS}
\left\{
\begin{aligned}
&   \mu \Delta v - \na q \: = \:  \overline{\rho} \left( \pa_t v +  v \cdot \na v -  f
  \right), \quad  \div v  = 0, \quad x \in F(t), \\ 
& v\vert_{\pa S^i(t)}  = v^i(t) + \omega^i(t) (x - x^i(t))^\bot,  \quad
  v\vert_{\pa \Omega \cap \pa F(t)} =0.   
\end{aligned}
\right.
\end{equation}  
As the solids and the cavity do not touch  ($\delta(t)\ge\delta_0$), it is
standard to build a   solenoidal vector field $w(t,\cdot) \in
H^\infty(\Omega)$ such that  
$$ w(t, \cdot)\vert_{S^i(t)} \: = \:  v^i(t) + \omega^i(t) (x -
x^i(t))^\bot, \quad  
w(t, \cdot)\vert_{\pa \Omega} \: =  \: 0, $$
 with the estimate $ \| w(t, \cdot) \|_{H^s} \: \le \:  C_s \sum_i
 (|v^i(t)| + |\omega^i(t)|)$. Then,  
the function $u = v - w$ satisfies 
\begin{equation*}
\left\{
\begin{aligned}
& \mu \Delta u - \na q =  \overline{\rho} \left( \pa_t v +  v \cdot \na v - f \right)
  - \mu \Delta w, \quad  \div u  = 0, \quad x \in F(t), \\ 
& u\vert_{\pa F(t)} = 0. 
\end{aligned}
\right.
\end{equation*} 
As $F(t)$ is a $C^{1,\alpha}$ open domain, we can apply the estimates of
the previous section. If $q$ is normalized so that $\int_{F(t)}  q(t,\cdot)
= 0$, we have, by \eqref{controlBMO} and propositions \ref{StokesBMO},\ref{StokesSob}
\begin{equation*}
\begin{aligned}
  &  \: \| (\na u, q)(t) \|_{L^2(F(t))} \: +  \:  \| (\na u, q)(t)
   \|_{H^1(F_{\eps}(t))}  \: 
   +     \:    \|(\na u,  q)(t) \|_{\bmo(F(t))} \\
   + \: &  \: \|(\na u, q)(t) \|_{W^{s,\tau}(F(t))} \:
 \le \:  C \, \Bigl( \| \pa_t v(t) \|_{L^2(F(t))} \: + \: \| v \cdot \na
   v(t) \|_{L^2(F(t))} \\ 
   & \hspace{3.9cm} + \:   \: \| f(t) \|_{L^2(F(t))} \: + \:  \sum_i
   (|v^i(t)| + |\omega^i(t) |) \Bigr).  
  \end{aligned}
  \end{equation*}
 We remind that this bound holds for all $s,\tau$ such that $s < \alpha$, and
 $s \le 2/\tau$.  
Back to the original field $v$, and using the interpolation inequality
 \eqref{interpolation}, we get: for all finite $r$ 
\begin{equation} \label{apriori2}
\begin{aligned}
   & \: \| (\na v,q)(t) \|^2_{L^r(F(t))} \: +  \:  \| (\na v,q)(t)
  \|^2_{H^1(F_{\eps}(t))}  
\: + \:    \|(\na v,q)(t) \|^2_{\bmo(F(t))} \\
+ &  \: \| (\na v,q)(t) \|^2_{W^{s,\tau}(F(t))}    
 \:  \le \: C' \, \Bigl( \| \pa_t v(t) \|^2_{L^2(F(t))} \: + \: \| v \cdot
  \na v(t) \|^2_{L^2(F(t))} \\ 
 & \hspace{4.1cm}  + \: \: \| f(t) \|^2_{L^2(F(t))} \: + \: 
   \sum_i (|v^i(t)|^2 + |\omega^i(t) |^2) \Bigr). 
  \end{aligned}
  \end{equation}
 By a time integration of \eqref{apriori2} from $0$ to $t$, and  a linear
 combination with \eqref{apriori} and \eqref{apriori1}, we obtain 
\begin{equation}  \label{apriori3}
\begin{aligned}
& \int_0^t  
  \Bigl(  \| \pa_t v(s) \|^2_{L^2(\Omega)} \:  + \: \| \na v(s)
  \|^2_{L^r(\Omega)} \: +  \:  \|(\na v, q)(s) \|^2_{H^1(F_{\eps}(s))} \: + \: \|
  (\na v,q)(s)) \|^2_{\bmo(F(s))}  \\ 
 &  + \: \| (\na v,q)(s) \|^2_{W^{s,\tau}(F(s))}   \Bigr) \, ds  \:
    + \: \| \na v(t) \|^2_{L^2(\Omega)} \\
    &   \: \le \: C(T) \left( \| v_0 \|^2_{H^1(\Omega)} \: + \: \| f
  \|^2_{L^2((0,T) \times \Omega)} \: + \:  \int_0^t \| v \cdot \na v(s)
  \|^2_{L^2(\Omega)} 
   \, ds \right) 
\end{aligned}
\end{equation}
where $C(T)$ is an increasing function of $T$. To have a  closed estimate,
 it remains to handle the nonlinear term.  
 We split it into 
$$  \int_0^t \int_\Omega |v \cdot \na v|^2 \: = \:  \int_0^t \int_{F(s)} |v
 \cdot \na v(s) |^2 \, ds   \: + \:  \sum_i \int_0^t \int_{S^i(s)} |v \cdot
 \na v(s) |^2 \, ds$$ 
The last term in the decomposition clearly satisfies
 $$ \sum_i \int_0^t \int_{S^i(s)} |v \cdot \na v(s) |^2 \, ds \: \le \: C
 \, \sum_i  \int_0^t |v_i(s)|^4 + |\omega^i(s)|^4 \, ds \: \le \: C'  $$  
 where $C'$ depends  on $\| v_0\|_{L^2}$ and  $\| f \|_{L^2((0,T)\times
   \Omega)}$.  The first term is bounded in the following way: 
 \begin{equation} \label{apriori4}
 \begin{aligned}
\int_0^t  \| v \cdot \na v(s)&  \|_{L^2(F(s))}^2 \,ds  \: \le \: \int_0^t
\| v(s) \|_{L^4(F(s))}^2 \| \na v(s) \|_{L^4(F(s))}^2 \, ds \\ 
 & \le \: C \,  \int_0^t \| v(s) \|_{L^2(F(s))} \| v(s) \|_{H^1(F(s))} \|
 \na v(s) \|_{L^2(F(s))} \| \na v(s) \|_{\bmo(F(s))} \, ds \\ 
  & \le \:  C \, \| v \|_{L^\infty(0,T; \, L^2(\Omega))} \int_0^t \| v(s)
  \|_{H^1(\Omega)}^2 \,  
  \| \na v(s) \|_{\bmo(F(s))} \, ds \\
  & \le \: C' \: +  \: \eta \int_0^t \| \na v(s) \|_{\bmo(F(s))}^2  \, ds
  \: +  \: C_\eta \, \int_0^t \| \na v(s) \|_{L^2(\Omega)}^4  
 \end{aligned}
 \end{equation}
 where $\eta$ is arbitrary  and $C'$ is  an increasing function of $\|
 v_0 \|_{L^2(\Omega)}$ and  $\| f \|_{L^2((0,T) \times \Omega)}$. 
 Note that the second line is deduced from the use of Gagliardo-Nirenberg
 inequality and the interpolation inequality \eqref{interpolation}. 
Choosing $\eta$ small enough, \eqref{apriori4} and \eqref{apriori3} imply that 
 $$   \| \na v(t) \|_{L^2(\Omega)} \: \le \: C + \int_0^t \| \na v(s) \|^2_{L^2(\Omega)} \, 
\| \na v(s) \|^2_{L^2(\Omega)} \, ds. $$
using Gronwall lemma, and the fact that $\int_0^T  \| \na v(s)
\|^2_{H^1(\Omega)}$ is bounded through \eqref{apriori}, we obtain  
$$ \| \na v \|_{L^\infty(0,T,L^2(\Omega))} \: \le \: C $$
where $C$ is an increasing function of $T$, $\| v_0\|_{H^1(\Omega)}$ and $
\| f \|_{L^2((0,T) \times \Omega)}$. Using this bound in  \eqref{apriori3},
we finally  obtain: \begin{equation} 
  \begin{aligned}
\int_0^T & 
 \Bigl( \| \pa_t v(t) \|^2_{L^2(\Omega)} \:  + \: \| \na v(t)
 \|^2_{L^r(\Omega)} \: +  \:  \| (\na v, q)(t) \|^2_{H^1(F_{\eps}(t))} \: +
 \: \| (\na v,q)(t) \|^2_{\bmo(F(t))} \\ 
 & + \: \| (\na v,q)(t) \|^2_{W^{s,\tau}(F(t))}  \Bigr) \, dt  
   \: + \: \| \na v \|_{L^\infty(0,T,L^2(\Omega))}   \: \le \: C\Bigl(T,\|
   v_0\|_{H^1(\Omega)},  \| f \|_{L^2(0,T \times \Omega)} \Bigr).  
  \end{aligned}
\end{equation}

\medskip
These a priori estimates are as usual  the key element in the construction
of strong solutions, as it provides compactness for a sequence of
approximate solutions. In the case of $C^{1,1}$ boundaries,  the issue of
building such approximate solutions and passing to the limit  has been
adressed  in B. Desjardins and M. Esteban, as well as in many other
studies. As it adapts straightforwardly to our case, we do not give further
detail and refer to these papers.  

\medskip
Let us stress that the $W^{s,\tau}$ regularity of $(\na v,q)$  allows to
define the stress tensor at the boundary  $\: \left( \pa_n v \, - \,  q\,
  n\right)\vert_{\pa F(t)}$. Indeed, taking indices $s,\tau$ such that $ \tau \,
s > 1$   (together with the requirements   $s < \alpha, \quad \tau \, s \le
2$), one can define the traces of $\na  v$ and $q$ as elements of
$W^{s-1/\tau,\tau}(\pa F(t))$ for almost all $t$.   
 Note also that  the regularity properties 
 $$v \in L^2(0,T; \, W^{1,4}(\Omega)), \quad \pa_t v \in L^2(0,T; L^2(\Omega))$$
 of a strong solution $v$ are enough to ensure that the right-hand side in
 \eqref{NSS} belongs to $L^2((0,T) \times \Omega)$. If $\delta(t) \ge
 \delta_0$ for all $t < T$, this automatically implies the
 $L^2(H^2_{loc})$, $L^2(\bmo)$ and $L^2(W^{s,\tau})$  bounds on  
$(\na v,q)$ restricted to  the fluid domain. This shows the last statement
of the proposition, and concludes the existence part.

\bigskip
{\em We now  turn to the uniqueness of strong solutions}. Our result is 
\begin{proposition} {\bf (Uniqueness of strong solutions)} 

\smallskip 
Let $\delta_0 > 0$, $\displaystyle 
\: v_0 \in H^1_0(\Omega), \quad \rho^i_0 \, D(v_0) = 0, \: \forall \, i, 
 \quad  \: f \in L^2(0,T;  W^{1,\infty}(\Omega)),  \:  \forall \, T >0.$ 
 Assume \eqref{nocontact}, and
$$\displaystyle \pa \Omega \in C^{1, \alpha}, \quad \pa S^i_0 \in C^{1, \alpha},
\quad \forall \, i, \quad 0 < \alpha \le 1.$$
There is at most one strong solution on $(0,T)$ such that $\delta(t) >
\delta_0 $ for all $t < T$.  
\end{proposition}

\medskip
{\bf Proof of the proposition.}  We follow closely the work of T. Takahashi
related to  $C^{1,1}$ boundaries. We focus on changes due to our not so
regular $C^{1,\alpha}$ domains.  As in \cite{Takahashi03bis}, we just consider
the case $N=1$, $f=0$,   that is one solid $S(t)$ immersed in the cavity
$\Omega$, without forcing. To lighten notations, we also assume that the
density $\rho=1$ in the solid and the fluid domains. Minor changes allow to
handle the general case.

\medskip
{\em Step 1:  Lagrangian coordinates.} The first step in  the analysis of
uniqueness for this free surface problem is to get back to a fixed domain,
by a  change of variables of lagrangian type. Let $v_0 \in H^1(\Omega)$ and
$S(0)$ the initial velocity field and solid position. We will denote by
$h(t)$ the position of the center of mass of the solid at time $t$. We can
always assume that $h(0) = 0$. Let $(v,q)$ a strong solution on $(0,T)$
such that $\delta(t) > \delta_0$ for all $t < T$.  

\medskip
We consider the same change of variables as in \cite[paragraph 4.1,
p1504]{Takahashi03bis}: as $\delta(t) > \delta_0$, a solenoidal velocity field
$\Lambda(t,x)$ is defined  such that   
\begin{align*} 
 \Lambda(t, x)    &   = 0,  & \mbox{ for } \: x  \: \mbox{ in  an  } \:
 \delta_0/4 \: \mbox{ neighborhood of } \pa \Omega,  \\ 
\Lambda(t,x) \:   & =    \dot{h}(t) + \omega(t) (x - h(t))^\bot,  & \mbox{
  for } \: x  \: \mbox{ in  an  } \: \delta_0/4 \: \mbox{ neighborhood of }
S(t). 
\end{align*}
 Then, one considers the flow 
$$ X(t, \cdot ) : \Omega \rightarrow \Omega,  \quad \frac{\pa}{\pa t}
X(y,t) = \Lambda(t, X(t,y)),  \quad X(0,y) \: = \: y. $$  
which maps $S(0)$ to $S(t)$ and $F(0)$ to $F(t)$. More precisely, in a
neighborhood of $S(0)$,  
$$X(t, y) = {h}(t) \: +  \: R_{\theta(t)} \, y, \quad \theta(t) =
\int_0^t \omega(s) \, ds, \quad R_\theta \: = \: \left( \begin{smallmatrix}
    \cos\theta & -\sin\theta \\ \sin \theta & \cos \theta\end{smallmatrix}
\right). $$ 
and near $\pa \Omega$, $X(t,y) = y$. Note that, as $(v,q)$ is a strong
solution, $h, \theta \in H^2(0,T)$. The mapping $X$ inherits the regularity
estimate 
$$ \| \pa_t^{i} X(t, \cdot) \|_{H^s(\Omega)} \: \le \: C_s \, ( |
h^{(i)}(t)|   + | \theta^{(i)}(t)| ), \quad \forall \, i=0,1,2, \quad
\forall \, s \in \N. $$ 
We then introduce the new functions 
$$ u(t,y) \: :=  \: (\na Y)^t(t, X(t,y)) \, v(t,X(t,y)), \quad p(t,y) \: :=
\: q(t, X(t,y)). $$ 
where $Y := X^{-1}$ denotes the inverse of $X$ with respect to the space
variable, and as usual $(\na Y)_{ij} = \pa_{x_i} Y_j$.  

\medskip
Following \cite[paragraph 4.2, p1507]{Takahashi03bis}, equations
\eqref{NS}-\eqref{contin2} turn into  
\begin{equation*}
\left\{
\begin{aligned}
& \pa_t u  +  M u  +  N u  -  \mu L u  +  G  p  =   0, \quad y  \in F(0), \\
& \div u \: =  \: 0, \quad y \in F(0), \\
& u(y,t) \: = \:  R_{-\theta(t)} \dot{h}(t) \: + \: \omega(t) y^\bot, \quad
y \in \overline{S(0)}, \\ 
& m \ddot{h}(t) \:  =  \: R_{\theta(t)} \, \int_{\pa S(0)} (\mu \na u - p) n \, dy, \\
& J \dot{\omega}(t) =  \int_{\pa S(0)}   (\mu \na u - p) n \cdot y^\bot \, dy,
\end{aligned}
\right.
\end{equation*}
plus the initial condition 
$$ u\vert_{t=0}(y) \: = \: u_0(t,y) \: := \: v_0(t, X(t,y)). $$
  We refer to \cite{Takahashi03bis} for the exact expression of the various
  operators. In short, $(\pa _t  + M) u $ corresponds to the original time
  derivative $\pa_t v$,  $\: N u$ corresponds to $v \cdot \na v$, $\: L u$
  corresponds to $\Delta v$, and $G p$ corresponds to $\na p$. {\em An
    important point is that 
  \begin{equation} \label{operators}
 N u = u \cdot \na u, \quad L u = \Delta u, \quad G p = \na p \quad \mbox{
   near } \: \pa \Omega \: \mbox{ and } \: S(0). 
\end{equation}}
Indeed, we have $X(t,y) = y$ near $\pa \Omega$, so that the change of
variables is trivial near the boundary of the cavity. Similarly, $\na
X(t,y) = R_{\theta(t)}$ near $S(0)$. As Navier-Stokes equations are
rotationally invariant, we get \eqref{operators}. 

\medskip
{\em Step 2: Stokes like formulation.} The operators above involve the flow
$X(t, \cdot)$, which depends on  the solution $u$ itself: hence,  they are
nonlinear. But  $X(0,y)=y$, which means that for small time, nonlinearities
are expected to be small. We shall therefore treat these nonlinear
perturbations as source terms. We introduce  
$$ H(t) \: := \: \int_0^t R_{-\theta(s)} \dot{h}(s) \, ds $$
and write the system as: 
\begin{equation} \label{stokeslike}
\left\{
\begin{aligned}
& \pa_t u   -  \mu \Delta u  +  \na  p  =   f  -  M u  -  u \cdot \na u,
\quad y  \in F(0), \\ 
& \div u \: =  \: 0, \quad y \in F(0), \\
& u(y,t) \: = \:   \dot{H}(t) \: + \: \omega(t) y^\bot, \quad y \in \overline{S(0)}, \\
& m \ddot{H}(t) \:  =  \:  \, \int_{\pa S(0)} (\mu \na u - p) n \, dy \ + \: w(t), \\
& J \dot{\omega}(t) =  \int_{\pa S(0)}   (\mu \na u - p) n \cdot y^\bot \, dy.
\end{aligned}
\right.
\end{equation}
where 
$$ f \: := \:   - \: (N u - u \cdot \na u) + \mu (L-\Delta) u - (G -\na) p, $$
and 
$$ w(t) \: = \: m \, \omega(t) R_{\theta(t)} \dot{h}(t)^\bot. $$

\medskip
{\em Step 3:  Uniqueness}. The uniqueness of the strong solution will be established thanks to the formulation \eqref{stokeslike}. Let $(v^1, q^1)$, $(v^2, q^2)$ two strong solutions  on $(0,T)$, $T > 0$, corresponding to the same initial velocity field $v_0 \in H^1(\Omega)$ and same initial configuration $S(0)$, $F(0) = \Omega\setminus \overline{S(0)}$. We remind that for the sake of brevity,  we consider the force-free case.  We assume that the boundaries $\pa \Omega$ and $\pa S(0)$ are $C^{1,\alpha}$, and that 
$$ \delta^1(t) > \delta_0, \quad \delta^2(t) > \delta_0, \quad \forall \, t \in [0,T).$$
We can associate to $v^i$ the change of variable $X^i$, the new functions $u^i$, $p^i$ and so on.
 We shall prove that $u^1 = u^2$ on $(0,T) \times \Omega$. The differences 
 $$u \: :=  \;  u^1 - u^2, \quad  H \:  := \:  H^1 - H^2, \quad \omega \:
 :=  \: \omega^1 - \omega^2$$ 
  satisfy with obvious notations
 \begin{equation} \label{stokeslike2}
\left\{
\begin{aligned}
& \pa_t u   -  \mu \Delta u  +  \na  p  =   f^1 - f^2 \: + \: M^2 u^2 - M^1
u^1  \: + u^2\cdot \na u^2 - u^1 \cdot \na u^1 , \quad y  \in F(0), \\ 
& \div u \: =  \: 0, \quad y \in F(0), \\ 
& u(y,t) \: = \:   \dot{H}(t) \: + \: \omega(t) y^\bot, \quad y \in \overline{S(0)}, \\
& m \dot{H}(t) \:  =  \:  \, \int_{\pa S(0)} (\mu \na u - p) n \, dy \ + \: w^1(t) - w^2(t), \\
& J \dot{\omega}(t) =  \int_{\pa S(0)}   (\mu \na u - p) n \cdot y^\bot \, dy.
\end{aligned}
\right.
\end{equation}
with initial condition $u\vert_{t=0} = 0$. 
Now, we can perform the exact same estimates  as those  performed earlier
to show existence of strong solutions. 
In particular, we get (see estimate \eqref{apriori3})
\begin{equation}  \label{totalbound}
\begin{aligned}
\int_0^T & 
 \Bigl( \| \pa_t u(t) \|^2_{L^2(\Omega)} \:  + \: \| \na u(t)
 \|^2_{L^r(\Omega)} \: +  \:  \| \na u(t) \|^2_{H^1(F_{\eps}(0))} \: + \:
 \| \na u(t) \|^2_{\bmo(F(0))}  \\ 
 & + \: \| p(t) \|^2_{H^1(F_{\eps}(0))} \Bigr) \, dt \: + \:  \| \na u
 \|_{L^\infty(0,T; \, L^2(\Omega))}   \\ 
 &  \le \: \textrm{RHS} \: :=  \: C(T) \Bigl( \| f^1 - f^2 \|^2_{L^2((0,T)
   \times \Omega)} \: + \:   \| u^1 \cdot \na u^1 - u^2 \cdot \na u^2
 \|^2_{L^2((0,T) \times \Omega)} \,  \\ 
 & + \: \| M^1 u^1 - M^2 u^2 \|^2_{L^2((0,T) \times \Omega)} \: + \: \| w^1 
 - w^2 \|_{L^2(0,T)} \Bigr)  
  \end{aligned}
\end{equation}
where $\eps$ is any constant lower than $\delta_0/4$, and $C(T)$ an
increasing function of $T$. As usual, the pressure $p$ is normalized so
that $\int_{F(0)} p =0$. It remains to estimate the right-hand side. 

\medskip
By the remark \eqref{operators}, $f^1 - f^2$ has a support $F_{\eps}$ which is
compact in $F(0)$. 
The $L^2(H^1)$ bound on $(\na u, p)$, which was true up to the boundary for
$C^{1,1}$ domains,   holds  in $F_{\eps}$. Moreover, the $L^2((0,T) \times
\Omega)$ estimate on $\pa_t u$  and the $L^\infty(0,T; \,  L^2(\Omega))$
estimate on $\na u$  also  hold.  Hence, the same bounds as those derived
in \cite[Corollary 6.16, p1523]{Takahashi03bis} apply:  
\begin{align*}
 \| f^1 - f^2 \|_{L^2(0,T \times \Omega)} \: & \le \: {\cal C} \, T^{1/10}
 \: \Bigl( \| (\na u,p)\|_{L^2(0,T; H^1(F_{\eps}(0))}   
  \\
  & \hspace{1.5cm}   + \| \pa_t u \|_{L^2((0,T) \times \Omega)}  \: + \: \|
  \na u \|_{L^\infty(0,T; L^2(\Omega))} \Bigr)  
 \end{align*}
where ${\cal C}$ denotes here and in the sequel an increasing function of
$T$ and $\|v_0 \|_{H^1(\Omega)}$.  

\medskip
As the $L^2$ bound on $\pa_t u$ is still available, we deduce  as in
\cite[corollary 6.16]{Takahashi03bis}, that 
$$ \| w^1 - w^2 \|_{L^2(H^2)} \: \le \: {\cal C} \, T^{1/2} \, \| \pa_t u
\|_{L^2((0,T) \times \Omega)}.  $$  

\medskip
We remind that $M u \: = \; \omega(t) \, u^\bot \: + \: \pa_t Y \cdot \na
u, \:$. Therefore, 
\begin{align*}
&  \| M^1 u^1  -  M^2 u^2 \|_{L^2((0,T) \times \Omega)}   \:  \le \: \| M^1
u \|_{L^2((0,T) \times \Omega)} \:  + \:  \| M u^2 \|_{L^2((0,T) \times
  \Omega)} \\ 
 & \: \le \: C \biggl(  \| (\omega^1, \dot{h}^1) \|_{L^2(0,T)} \| u \|_{L^2(0,T;
   H^1)} \: + \:  \| (\omega, \dot{h}) \|_{L^2(0,T)} \| u^2 \|_{L^2(0,T; H^1)}
 \biggr)\\ 
 &  \: \le \: {\cal C}  \, T^{1/2} \, \| \na u \|_{L^\infty(0,T; H^1(\Omega))}.  
 \end{align*}
 
 \medskip
 Finally, we must control the quadratic term 
 $$u^1 \cdot \na u^1 - u^2 \cdot \na u^2 \: = \: u \cdot \na u^1 + u^2 \cdot \na u. $$
 Like in previous computations, we get 
 \begin{align*}
  \|  u \cdot \na u^1 \|^2 _{L^2((0,T) \times F(0))}  & \le C \|  u
  \|^2_{L^\infty(0,T; \, L^4(\Omega))}  \| \na u^1
  \|_{L^\infty(L^2(\Omega))} \int_0^T \| \na u^1(t) \|_{\bmo(F(0))}      \\
  & \le \:  {\cal C} \: \sqrt{T} \,  \| \na u \|_{L^\infty(0,T; H^1(\Omega))}^2, 
  \end{align*} 
 using the $L^2(\bmo)$ bound on $\na u^1$.  Similarly, 
 \begin{align*}
  \|  u^2 \cdot \na u \|^2 _{L^2((0,T) \times F(0))}  & \le C \|  u^2
  \|^2_{L^\infty(L^4(\Omega))}  \| \na u \|_{L^\infty(L^2(\Omega))}
  \int_0^T \| \na u(t) \|_{\bmo(F(0))}      \\
  & \le \: {\cal C} \, \sqrt{T}  \, \left( \| u
    \|_{L^\infty(H^1(\Omega))}^2 \: + \: \| \na u \|_{L^2(\bmo(F(0)))}
  \right).  
  \end{align*} 
 The $L^2$ bound in $S(0)$ is straightforward, and we end up with 
 \begin{equation*}
 \| u^1 \cdot \na u^1 - u^2 \cdot\na u^2 \|_{L^2(0,T \times \Omega)} \: \le
 \: {\cal C} \, T^{1/4} \,  \left( \|  u \|_{L^\infty(H^1(\Omega))}^2 \: +
   \: \| \na u \|_{L^2(\bmo(F(0)))} \right) 
 \end{equation*}
  
 \medskip
Eventually, these inequalities lead to
\begin{align*}
 \textrm{RHS} \: &  \le \: {\cal C} \, T^{1/10} \, \Bigl(  \| u
 \|_{L^\infty(H^1(\Omega))} \: + \:  \| \pa_t u \|_{L^2(L^2(\Omega))} \\ 
 &  \hspace{1.5cm} + \: \| \na u \|_{L^2(\bmo(F(0)))} \: + \:  \| (\na u,
 p) \|_{L^2(H^1(F_{\eps}))} \Bigr)  
 \end{align*}
with ${\cal C}$ an increasing function of $T$ and $\| v_0
\|_{H^1(\Omega)}$.  By reporting this bound in  \eqref{totalbound}, we
deduce that there exists a small $T_0$ such that  $u^1 = u^2$ on $[0,
T_0]$. Moreover,   
$T_0$ depends only   on $\| v_0 \|_{H^1(\Omega)}$ (decreasing as $\| v_0
\|_{H^1(\Omega)}$ increases).  

\medskip
{\em Global uniqueness follows}. Indeed, we know that $u^1$ and $u^2$ are
in $L^\infty([0,T], H^1(\Omega))$ for all times $T$ such that  $ \forall \,
t
\le T, \; \delta^1(t), \delta^2(t) \ge \delta_0$. Hence, up to consider a
smaller $T_0$, we can apply the above local uniqueness result on $[T_0,
2T_0]$, then $[2T_0, 3 T_0]$ and so on up to reach time $T$. This concludes
the uniqueness part.  

\bigskip
{\em So far, we have shown existence and uniqueness of strong solutions at
  least up to the first collision}. Theorem \ref{theo1} asserts more,
namely that no strong solution can exist beyond the first collision
time. This result can be deduced from \cite[Theorems 3.1 and 3.2b]{Starovoitov03}. 
Indeed, V. Starovoitov has shown the following: suppose
that    two $C^{1,\alpha}$ solids $S^{i_1}(t)$  and $S^{i_2}(t)$, {\it
  resp.}   a $C^{1,\alpha}$ solid $S^i(t)$ and the $C^{1,\alpha}$ cavity
$\Omega$,  collide for the first time at   $t = T$. Denote   for $t \le T$,  
$$h_{i_1,i_2}(t) \: := \: \mbox{dist}\left(S^{i_1}(t), S^{i_2}(t)\right)) \quad
\mbox{{\it resp.} } \: h_{i}(t) \: := \: \mbox{dist}\left(S^{i}(t),
  \Omega\right) $$ 
 and assume that 
 \begin{equation} \label{regularuptocoll}
 u \in L^\infty(0,T; \, L^2(\Omega)) \cap L^1(0,T; \, W^{1,p}(\Omega)).
 \end{equation}
Then, $h_i,h_{i_1,i_2}$ are lipschitz continuous on $[0,T]$ and, for example, 
$$ \left|\frac{dh_i}{dt}(t)\right| \: \le \: C \, h_i(t)^\beta \, \| u(t)
\|_{W^{1,p}(\Omega)}, \quad \beta = 2 - \frac{1}{1+\alpha} \, \frac{p+1}{p}
\: - \: \frac{1}{p}, $$ 
for almost all $t \le T$. 

\medskip
 In particular, if $T$ is the first collision time, and the strong solution
 exists beyond it, the regularity assumption \eqref{regularuptocoll}  is
 satisfied for arbitrary $p$. Taking $p$ large enough, one can assume that
 $\beta \ge 1$. We note also that, by hypothesis \eqref{nocontact}, $h(0)
 \neq 0$. Then, by integration of  the previous differential inequality, we
 obtain $h(T) \neq 0$, which yields a contradiction. 

We emphasize that for strong solutions, 
\eqref{regularuptocoll} holds {\it a priori} only for $T<T_*.$
Thus, the argument of Starovoitov does not allow to conclude on the occurence
of collision. In the next section, we will exhibit configurations for which
collision occurs. For such examples, we have:
$$
\int_0^{T_*} \|u\|_{W^{1,p}(\Omega)} = \infty. 
$$

\section{The collision/no-collision result} \label{secconst}

%%%%%%%%%%%%%%%%%%%%%%%%%%%%%%%%%%%
%%
%%				Section Collision solide non regulier
%%
%%%%%%%%%%%%%%%%%%%%%%%%%%%%%%%%%%%

This section is devoted to the proof of Theorem \ref{theo2}. We consider the simplified configuration described at the end of the introduction, assumptions 1-6. In this framework, the position of the solid is characterized by 
$ h(t) \: := \: \textrm{dist} \left( (0,x_-(t)), \, \pa \Omega \right).$
Later on, it will be convenient to use a parametrization by $h$, {\it i.e.} the translated domains 
$$ S_h \: := \: S(0) \: + \: (h - h(0)) e_2, \quad h \in \R. $$
Of course, $S_{h(t)} = S(t)$. By assumption 5, the boundary of $S(t)$ is
$C^{1,1}$ near its "upper tip" $(0,x_+(t))$, so that contact is impossible
at this point, {\it cf.} \cite{Hillairet07}. By assumption 6, gravity pushes
$S(t)$ downwards. So, we can even assume that   
\begin{equation}\label{nocontact2}
 \inf_{t \in (0,T_*)} \, \textrm{dist} \left( (0,x_+(t)),
   \pa \Omega \right) \: > 0.   
 \end{equation} 

\medskip
Thus, collision can occur if and only if $\lim_{t\rightarrow T_*} h(t) =
0$. We will show that it is equivalent to $\alpha < 1/2$. The proof is
based on the use of a quasistationary velocity field $w$ and
quasistationary pressure field $q$.  By quasistationary, we mean that for
all $t < T_*$,  
$$ w(t,x) \: = \: w_{h(t)}(x),  \quad q(t,x) \: = \: q_{h(t)}(x)  $$
for some stationary fields $w_h(\cdot)$, $p_h(\cdot)$ defined on $\Omega$
and parametrized by $h > 0$. Moreover, they will satisfy  
\begin{equation} \label{regularityw}
  w \in C^1([0,T_*); \, H^1(\Omega)), \quad  \Delta w(t, \cdot) \in
  L^p(F(t)), \:   p \mbox{ small enough},  \: t \in (0,T_*),  
\end{equation}
\begin{equation} \label{propertyw}
\div w = 0 \mbox{ in } \Omega,  \quad w\vert_{S(t)} = e_2, \quad w\vert_{\pa \Omega} = 0. 
\end{equation}
as well as 
\begin{equation} \label{regularityq}
q \in C^1([0,T_*), L^2(\Omega)), \quad  \na q(t, \cdot) \in L^p(F(t)), \:
p \mbox{ small enough},  \: t \in (0,T_*). 
\end{equation}
In particular, we can use $w$ as a test function in the variational formulation to get 
\begin{align*}
&  \int_0^t \int_\Omega \left( \rho v \cdot \pa_t w + \rho v \otimes v :
  D(w) - 2 \mu D(v) : D(w)  -\rho g e_2 \cdot w \right) \\ 
 & = \int_\Omega \rho(t) v(t) \cdot w(t) - \int_\Omega \rho(0) v_0 \cdot w(0) 
\end{align*}
Note that by \eqref{propertyw}
\begin{align*}
\int_{\Omega} \rho g e_2 \cdot w \: & = \: \int_{S(t)} \rho g e_2 \cdot e_2 \: + \: \int_{F(t)} \rho g \na (x \mapsto x_2) \cdot w \\
& = \: \rho_S \, g \, |S(0)| \: + \: \rho_F \, g \, \int_{\pa F(t)} x_2 \, e_2 \cdot n \:  = \:  \left( \rho_S - \rho_F \right) \, g | S(0)| 
\end{align*}
where $\rho_S := \rho\vert_{S(t)}$, $\: \rho_F := \rho\vert_{F(t)}$.  We also write 
\begin{align*}
 2 \mu \int_{\Omega} D(v)(t) : D(w)(t) \: & = \: \dot{h}(t) \int_{\pa F(t)}  \left(\mu \frac{\pa w}{\pa n} - q n \right) \cdot e_2 \: - \: \int_{F(t)} (\mu \Delta w - \na q)   v  \\
  & := \: \dot{h}(t) \, n(h) \: - \:  \int_{F(t)} (\mu \Delta w - \na q)  \cdot  v.
  \end{align*}
Thus, the variational formulation yields
\begin{equation} \label{eqN(h)} 
 N(h(t)) \: + \:  \left( \rho_S - \rho_F \right) \, g | S(0)| \, t \: = \: R(t) 
\end{equation}
where $N$ is the antiderivative of $n$ that vanishes at $h(0)$, and the remainder is 
\begin{align*}
 R(t)\: & := \:  \int_0^t \int_\Omega \left( \rho v \cdot \pa_t w + \rho v \otimes v : D(w) \right) 
 \: + \: \int_\Omega \rho(0) v_0 \cdot w(0) - \int_\Omega \rho(t) v(t) \cdot w(t) \\
 & \:  + \: \int_0^t  \int_{F(s)} \Delta w(s, \cdot) - \na q(s, \cdot)   \cdot v(s, \cdot). 
 \end{align*}
 Theorem \ref{theo2} will be deduced from the following proposition:
 \begin{proposition} \label{propn(h)}
 One can find $w_h : \Omega \mapsto \R^2$, $\: q_h : \Omega \mapsto \R$,
 such that $w,q$ satisfy \eqref{regularityw},  \eqref{propertyw}, \eqref{regularityq},
 and such that 
 \begin{description} 
 \item[i)] For $h > 0$ small enough 
 \begin{equation} \label{estimn(h)}
 -c \:  \le \:  n(h) \: \le \: C \, h^{-\beta}, \quad \beta =
 \frac{3\alpha}{1+\alpha}, \quad c, C > 0.  
 \end{equation}
 \item[ii)] For all $t < T_*$,
 \begin{equation} \label{estimrm}
 |R(t)| \le C(\| u_0 \|_{L^2}) \left( 1 + \sqrt{t} \right). 
 \end{equation}
 \end{description}
 \end{proposition}
 Before tackling the proof of this proposition, let us show how it implies Theorem \ref{theo2}. 
 
 \medskip
 {\em If $\alpha \ge 1/2$}, then $\beta \ge 1$. We get from \eqref{eqN(h)}
 and point ii) of the proposition: 
 $$ N(h(t)) \: \ge \:  (\rho_F -
  \rho_S) |S(0)| t \: - C \, (1+\sqrt{t}) $$
 By point i), we also get for $h$ small enough 
 $$ N(h) \: \le \:   - C | \ln(h) | $$
 In fact, one can take $h^{1 - \beta}$ instead of $|\ln(h)|$ when $\beta >
 1$, {\it i.e.} $\alpha > 1/2$. Combining those inequalities, we deduce  
 $$
 C | \ln h(t) | \; \le \:  (\rho_S -
  \rho_F) |S(0)| t \: + C (1 + \sqrt{t}) \: <  \: + \infty, \quad
 \forall \, t < T_* 
 $$ 
 which means that $h$ does not go to zero in finite time. Hence,
 $T_*=+\infty$ and there is no collision  
 
 \medskip
  {\em If $\alpha < 1/2$}, then $\beta < 1$, and $n \in L^1$. Thus, $N$ is
  continuous. As $h(t)$ is bounded, we deduce from \eqref{eqN(h)}: $\forall
  \, t < T_* \le +\infty$,  
  $$  -\infty \: < \: \inf_{t \in (0,T_*)} N(h(t)) \: \le \:  \:  (\rho_F -
  \rho_S) |S(0)| t \: + \: C (1 + \sqrt{t}). 
  $$  
  If $T_*=+\infty$, and $\rho_S > \rho_F$, one can let $t \rightarrow
  +\infty$ in the previous inequality. As the r.h.s. goes to $-\infty$ in
  this limit, it yields a contradiction. Thus, $T_* < +\infty$. This ends
  the proof.  
   
\medskip
The rest of the paper will be devoted to the proof of Proposition \ref{propn(h)}.

\subsection{Construction of the test function}
We mimic the construction presented in article \cite{Hillairet07} for
$C^{1,1}$ boundaries. We want a function $w_h(x)$ such that  
\begin{equation} \label{propertywh}
 \div w_h \: = \: 0 \mbox{ in } \Omega, \quad w_h\vert_{S_h} \: = \: e_2,
 \quad w_h\vert_{\pa \Omega} \: = \: 0.  
 \end{equation}
 We always consider $0 < h < h_M := \sup_{0<t<T_*} h(t)$, as no other value
 of $h$ is involved in our problem.  
 
\begin{figure} \label{geom3}
\begin{center}
\includegraphics[height = 5cm, width=5.5cm]{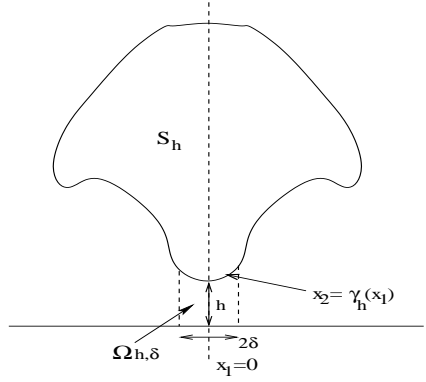}
\end{center}
\caption{Geometry of the possible contact zone}
\end{figure}

 \medskip
 By a change of coordinates, we can assume $(0,0) \in \pa \Omega$, {\it
   i.e.} $x_-(t) = h(t)$. By assumption 5, there exists $\delta > 0$, such
 that  
 $$ \forall \, x \in \pa S_h \cap D((0,h), 2\delta), \quad x_2 = \gamma_h(x_1) \: :=  \: 
  h + |x_1|^{1+\alpha}, $$
 where as usual $D(x,r)$ is the disk of center $x$ and radius
 $r$. Moreover, by assumption 3  and \eqref{nocontact2},  
 $$ \delta_{min}\: := \: \inf_{0 , h < h_M} \textrm{dist}\left( \pa S_h \cap
   D((0,h), \delta)^c, \, \pa \Omega \right) \: >  \: 0.  $$ 
 To describe $w_h$ away from the origin, we introduce a smooth function 
 $\varphi = \varphi(x)$, $x \in \R^2$ such that 
 \begin{equation*}
 \varphi = 1  \: \mbox{ in a $\delta_{min}/2$-neighborhood of } S_{h(0)}, \quad
 \varphi = 0  \: \mbox{ outside a $\delta_{min}$-neighborhood of } S_{h(0)}, 
 \end{equation*}
We introduce another smooth  function $\chi = \chi(x)$, $x \in \R^2$,  such that 
 \begin{equation*}
 \chi = 1  \: \mbox{ in } (-\delta, \delta)^2, \quad
 \chi = 0  \: \mbox{ outside } (-2\delta, 2\delta)^2.
 \end{equation*}
 Finally, we set $w_h = \na^{\bot}(x_1 \varphi_h)$, with 
 \begin{align*}
 \varphi_h & = 1  \quad \mbox{in } S_h, \\
 \varphi_h  & = (1-\chi(x)) \, \varphi(x_1, x_2 - h + h(0)) \:  +  \; \chi(x) \, 
 \frac{x^2_2}{\gamma_h(x_1)^2}  \left( 3 - \frac{2x_2}{\gamma_h(x_1)}\right)
 \quad \mbox{ in  } \Omega \setminus S_h.
 \end{align*}
 See figure \ref{geom3} to clarify the main  notations.  Note that $\varphi_h$ and
 therefore $w_h$ are regular up to $h=0$ outside  
 $$ \Omega_{h, \delta} \: := \: \Omega \cap \{ |x_1| < \delta \} \cap \{ x_2 < \gamma_h(x_1) \}. $$
 Singularities at $h=0$ correspond to  the second term in the definition of $\varphi_h$. 

 \medskip
 It is straightforward that $w_h$ satisfies \eqref{propertywh}. 
 As $\varphi_h$ involves the boundary function $\gamma_h$, the
 streamfunction $x_1 \varphi_h$  has  regularity $C^{2,\alpha}$ in the
 fluid domain. Moreover,  $ w_h$ is continuous across the solid boundary,
 so that it 
 belongs to $C^\infty((0, h_M); \,  W^{1,\infty}(\Omega))$. In the fluid domain, its
 most singular second order derivatives behave like $x_1^{\alpha-1}$. We
 deduce that  $w(t,x) = w_{h(t)}(x)$ satisfies \eqref{regularityw}.  
We postpone to  the appendix the proof of the  following estimates:
\begin{proposition} \label{estimatewh}
There exists $0 \: < \:  c \: < \: C$ such that 
\begin{equation} \label{estimwh1}
\begin{aligned}
   \| w_h \|_{L^2(\Omega)} \:  &  \le \: C, \\
 c \,\:  \le \:    h^{\frac{3\alpha}{2(1+\alpha)}}  \| \na w_h \|_{L^2(\Omega)} \: & \le \: 
C, \\
 \|\na w_h \|_{L^\infty\left(\Omega\setminus\Omega_{h,\delta}\right)} \: & \le \: C, 
\end{aligned}
\end{equation}
and 
\begin{equation} \label{estimwh2}
\begin{aligned}
 \sup_{x_1 \in (-\delta, \delta)} |\gamma_h(x_1)|^{3/2} \, \left( 
\int_0^{\gamma_h(x_1)} | \na w_h(x_1, x_2)|^2 dx_2  \right)^{1/2} \: & \le \: C, \\
 \int_{-\delta}^\delta \int_0^{\gamma_h(x_1)} \gamma_h(x_1)^2 \,  \: |\pa_h
 w_h(x)|^2   \, dx \; & \le \: C. 
\end{aligned}
\end{equation}
\end{proposition}
Besides these estimates on $w_h$, the control of $n(h)$ and $R(t)$ shall
involve quantities of the type  
$$ \int_{F(t)}(\mu \Delta w_h - \na q_h)   \tilde{w}  $$ 
where $\tilde{w} \in H^1_0(\Omega)$ is divergence free  and satisfies
$\tilde{w}\vert_{\pa S_h}  = e_2$. We prove in the appendix the following estimate 
\begin{proposition} \label{propqh}
There exists a pressure field $h \mapsto q_h \in C^\infty(0,h_M;
C^1(\overline{\Omega}))$ such that  for all divergence free $\tilde{w} \in
H^1_0(\Omega)$ satisfying  
 $\tilde{w}\vert_{\pa S_h}  = e_2$, 
$$ \left|  \int_{F(t)}(\mu \Delta w_h - \na q_h)   \tilde{w} \right| \: \le \:
C \, \| \tilde{w} \|_{H^1_0(\Omega)}. $$ 
\end{proposition}

\subsection{Proof of Proposition \ref{propn(h)}}
Thanks to the estimates of the previous section, we are able to control the
functions $n(h)$ and the remainder term $R(t)$.  

\medskip
By integration by parts,
$$ n(h) \: := \:  2\mu \int_{\Omega}D(w_h) : D(w_h) \: + \: \int_{F_h} (\mu \Delta w_h - \na q_h) w_h. $$
By estimate  \eqref{estimwh1}b, we have
$$ c   h^{-\frac{3\alpha}{2(1+\alpha)}} \,\:  \le \:   \| \na w_h \|_{L^2(\Omega)} \:  \le \: 
C   h^{-\frac{3\alpha}{2(1+\alpha)}} , $$
and Proposition \ref{propqh} leads to 
$$ \left|  \int_{F(t)}(\mu \Delta w_h - \na q_h)   w_h \right| \: \le \: C \,
\| \na w_h \|_{L^2(\Omega)} \; \le \: C_\eps \: + \: \eps  \| \na w_h
\|_{L^2(\Omega)}^2. $$ 
Combining these last two inequalities yields point i) of Proposition \ref{propn(h)}.

\medskip
To establish point ii), we need to control each term in the
remainder. Still using bound \eqref{estimwh1}, we have  
$$ \left|\int_{F_{h(t)}} (\mu \Delta w_{h(t)} - \na q_{h(t)}) \, v \right| \:
\le \: C \| v(t) \|_{H^1_0(\Omega)}.$$ 
Integration from $0$ to $t$ and Cauchy-Schwarz inequality lead to 
\begin{equation} \label{remainder1}
 \int_0^t \left|\int_{F(s)} (\mu \Delta w(s, \cdot)- \na q(s, \cdot)) \, v(s,
   \cdot) \right| \: \le \: C \, \| v \|_{L^2(0,t; H^1(\Omega))} \,
 \sqrt{t} \: \le \: C(\| u_0 \|_{L^2(\Omega)})  \, \sqrt{t}. 
 \end{equation}
 
\medskip 
We also get 
\begin{align}
\nonumber
 \left| \int_\Omega \rho(0) v_0 \cdot w(0) - \int_\Omega \rho(t) v(t) \cdot w(t) \right| \: & 
 \le \: C \, \| v \|_{L^\infty(0,t; L^2(\Omega))} \, \sup_{h \in (0,h_M)} 
\| w_h(t) \|_{L^2(\Omega)} \\
\label{remainder2}
 &  \le \: C(\| u_0 \|_{L^2(\Omega)}) 
\end{align} 

\medskip
To deal with the term involving $\pa_t w$, we shall use the following
general bound: for any $h \in (0,h_M)$ and any $(\rho,v) \in
L^{\infty}(\Omega) \times H^1_0(\Omega)$  we have, for any $\tilde{w} \in
H^1_0(\Omega):$ 
\begin{align*}
\left| \int_{\Omega} \rho v \cdot \tilde{w} \right| \leq C \|\rho\|_{L^\infty(\Omega)} \, 
 \|\nabla v\|_{L^2(\Omega)}  \biggl(  
& \|\tilde{w}\|_{L^2(\Omega \setminus \Omega_{h,\delta})} \\
& + \:  \biggl(\int_{-\delta}^\delta \int_{0}^{\gamma_h(x_1)} |\gamma_h(x_1)|^2 |\tilde{w}(x)|^2  
\text{d$x$}\biggr)^{1/2} \biggr),
\end{align*}
This is a simple consequence of Cauchy-Schwarz and Hardy inequalities, and
its proof is therefore left to the reader. Note that $\pa_t w = \dot{h}(t)
\pa_h w_{h(t)}$. The previous  formula  yields  
\begin{align}
\nonumber
&  \left| \int_0^t \int_\Omega  \rho v \cdot \pa_t w   \right|  \:  \le \:
C \sup_{[0,T_*)} |\dot{h}| \, \int_0^t \| \na v(s) \|_{L^2(\Omega)}
\biggl(   \|\pa_h w_{h(s)}\|_{L^2(\Omega \setminus \Omega_{h(s),\delta})}
\\ 
 \label{remainder3}
& + \:  \biggl(\int_{-\delta}^\delta \int_{0}^{\gamma_{h(s)}(x_1)}
|\gamma_{h(s)}(x_1)|^2 |\pa_h w_{h(s)}|^2   dx\biggr)^{1/2}
\biggr) \, ds \: \le \: C(\| u_0 \|_{L^2(\Omega)})    \sqrt{t}, 
\end{align}
where the last inequality involves \eqref{estimwh2}b.  
Finally, to deal with the nonlinear term, we use another general formula, namely: 
for any $h \in (0,h_M)$ and any $(\rho,v) \in L^{\infty}(\Omega) \times
H^1_0(\Omega)$  we have, for any $\tilde{w} \in H^1_0(\Omega): $ 
\begin{align*}
 \left| \int_{\Omega} \rho v\otimes v : D(\tilde{w}) \right|  \leq & C \|\rho\|_{L^\infty(\Omega)} 
\, \|\nabla v\|_{L^2(\Omega)}^2
\bigg( 
\|D(\tilde{w})\|_{L^\infty(\Omega \setminus \Omega_{h,\delta})}
\\ 
& + 
\sup_{x_1 \in (-\delta,\delta)} 
\biggl(
	|\gamma_h(x_1)|^{\frac{3}{2}}
	\biggl( 
		\int_{0}^{\gamma_h(x_1)} |\nabla \tilde{w}(x)|^2 \text{d$x_1$} 
	\biggr)^{\frac{1}{2}} 
\biggr)
\biggr).
\end{align*}
This formula follows from Cauchy-Schwarz inequality together with a refined
Poincar\'e's inequality. We refer to lemma 12 in \cite{Hillairet07} for all necessary
details. We infer from  this bound and \eqref{estimwh2} that 
\begin{equation} \label{remainder4}
   \int_0^t \int_\Omega  \rho v \otimes v : D(w)  \: \le \: C(\| u_0 \|_{L^2(\Omega)}). 
\end{equation}
Gathering \eqref{remainder1} to \eqref{remainder4} gives the bound on $R(t)$.

\section*{Acknowledgements}
The authors wish to thank Bertrand Maury for a useful discussion.  

\bibliographystyle{plain}
\bibliography{biblio}
\appendix
\section*{Appendix : Proofs of propositions \ref{estimatewh} and \ref{propqh}} \label{sec_w}
In this section, we estimate the rate of divergence of various Sobolev norms
of $w_h$ as $h$ goes to $0.$ As explained in section \ref{secconst}, $w_h$
is regular up to $h=0$ in $S_h$ and  $\Omega \setminus (S_h \cup
\Omega_{h,\delta}).$ Hence, there holds:
$$
\| \na w_h \|_{L^\infty(\Omega\setminus\Omega_{h,\delta})}   \leq C,
$$
and the rate of divergence of $w_h$ is the one of its
restriction to $\Omega_{h,\delta}$ \emph{i.e.} the one of $\nabla^{\bot}( x_1\varphi_h)$ 
where
$$
\varphi_h(x) =
\dfrac{x_2^2}{\gamma_h(x_1)}\left(3 - 2
  \dfrac{x_2}{\gamma_h(x_1)}\right),\quad
\forall \, x \in \Omega_{h,\delta}.
$$ 
Proposition \ref{estimatewh} is then a straightforward consequence of: 
\begin{lemma} \label{lem_size}
Given $(\alpha,p,q) \in (0,\infty)^3,$ the quantity:
$$
\int_{-\delta}^{\delta} \frac{|x_1|^p\text{d$x_1$}}{(h+|x_1|^{1+\alpha})^q}
$$
behaves like \vspace{.2cm}\\
\begin{tabular}{cll}
(i)  & $c h^{\frac{(p+1)}{1+\alpha} - q},$ & if $p+1 < q(1+\alpha),$ \\
(ii) & $c \ln(h),$                         & if $p+1 = q(1+\alpha),$ \\
(iii)& $c,$                                & if $p+1 > q(1+\alpha).$
\end{tabular}                                          \vspace{.2cm}\\
when $h$ goes to $0,$ with $c$ a constant depending only on $(\alpha,p,q).$
\end{lemma}
The proof of this lemma as well as the induced bounds on $w_h$ are direct
adaptation of \cite[Lemma 13]{Hillairet07}. 

\medskip
It remains to build the pressure field $q_h$ in order to prove proposition
\ref{propqh}. For simplicity, we assume now $\mu = 1.$ With the same notations as in section \ref{secconst}, we set:
$$
q_h(x) = \partial_{21} (x_1 \varphi_h(x))\:  + 12 \: \int_{0}^{x_1}  \frac{t}{\gamma_h(x_1)^2}
\text{d$t$}, \quad \forall \,  x \in \Omega.
$$
We stress that 
$$
q_h(x) = \partial_{21} (x_1 \varphi_h(x))\:  - \:
\int_{0}^{x_1} \partial_{222} (t \, \varphi_h(t,x_2))
\text{d$t$}, \quad \forall \, x \in \Omega_{h,\delta}.
$$
As for $w_h,$ this pressure field is smooth up to $h=0$ in the fluid domain
outside $\Omega_{h,\delta}.$ Consequently, the rate of divergence of
$\Delta w_h - \nabla q_h$ is the one of its restriction to this latter
domain. Standard computations lead to :
$$
\Delta w_h (x) - \nabla q_h(x) = 
\left(
\begin{array}{c}
-2 \partial_{112} (x_1 \,\varphi_h(x))\\
\partial_{111}(x_1 \varphi_h(x)) 
\end{array}
\right)
\quad \forall \, x \in \Omega_{h,\delta}.
$$
We recall that $\nabla^2 w_h \in L^p(\Omega_{h,\delta})$ for $p$
sufficiently small. As $H^1(\Omega_{h,\delta}) \subset 
L^r(\Omega_{h,\delta_0})$ for arbitrary $r < \infty,$ the integral to be
estimated in proposition \ref{propqh} is well-defined. 
Up to a truncation (which leaves aside a term that is regular with respect
to $h$), we can assume $\tilde{w} = 0$ in
$(\Omega \setminus S_h )\setminus \Omega_{h,\delta}$. 
 {\it A fortiori:}
$$
\int_{\Omega \setminus S_{h}} (\Delta w_h  - \nabla q_h)  \cdot \tilde{w}= 
\int_{\Omega_{h,\delta}} (\Delta w_h  - \nabla q_h)  \cdot \tilde{w}. 
$$ 
After an integration by parts, accounting for $\tilde{w}\vert_{\pa S_h}=
(0,\tilde{w}_2)$:
$$
\int_{\Omega_{h,\delta}} (\Delta w_h  - \nabla q_h)  \cdot \tilde{w} =
-\int_{\partial S_h} \partial_{11}(x_1 \, \varphi_h) \tilde{w}_2n_1 
\text{d$\sigma$} - \int_{\Omega_{h,\delta}}  \partial_{11} (  x_1 \,
\varphi_h ) 
(2 \partial_2 \tilde{w}_1 - \partial_1 \tilde{w}_2). 
$$
Thanks to Lemma \ref{lem_size}, one can check that $\pa_{11}  ( x_1 \,\varphi_h )$ is
bounded uniformly in $h$ in $L^2(\Omega_{h,\delta})$. 
 then shows  $\|\partial_{11} \varphi_h \|_{L^2(\Omega_{h,\delta})}$. Moreover, the boundary term reads
$$ 
\left|\int_{\partial S_h} \partial_{11}(x_1 \, \varphi_h) \, \tilde{w}_2 n_1
\text{d$\sigma$} \right| \:  \le \:  
\|\tilde{w}_2\|_{L^{\infty}(\partial \Omega_{h,\delta})} \int_0^\delta \left| \frac{6x_1(\gamma'_h(x_1))^2}{(\gamma_h(x_1))^2} \,
\frac{\gamma_h'(x_1)}{1+(\gamma'_h(x_1))^2} \right| \, dx_1 $$\
where $|\gamma'_h(x_1)| \leq c|x_1|^{\alpha}.$ So, this boundary term is again uniformly bounded by Lemma \ref{lem_size}. 
This ends the proof of proposition  \ref{propqh}.

\end{document}